\documentclass{article}
\usepackage[utf8]{inputenc}
\usepackage[small,compact]{titlesec}
\usepackage{amsmath,amsthm,amssymb,amsfonts,mathabx,mathrsfs,setspace}
\usepackage{graphicx,caption,epsfig,subfigure,epsfig,wrapfig,sidecap}
\usepackage{url,color,verbatim,algorithmicx}
\usepackage[sort,compress]{cite}
\usepackage{soul}
\usepackage{bm,float}
\usepackage[top=1in,bottom=1in,left=1in,right=1in]{geometry}
\usepackage{multirow, makecell}
\usepackage{cite}
\usepackage[ruled,boxed]{algorithm}
\usepackage[scaled]{helvet}
\usepackage[T1]{fontenc}
\usepackage[bookmarks=true, bookmarksnumbered=true, colorlinks=true,   pdfstartview=FitV,
linkcolor=blue, citecolor=blue, urlcolor=blue]{hyperref}
\usepackage{multirow, multicol,makecell, colortbl, hhline}
\usepackage{booktabs} 
\usepackage{enumitem} 
\usepackage{etoc}          

\numberwithin{equation}{section}

\newtheorem{theorem}{Theorem}[section]

\newtheorem{remark}{Remark}[section]

\newcommand{\mbf}[1]{\bm{#1}} 

\begin{document}

\title{Automated construction of effective potential via algorithmic implicit bias}

\author{Xingjie Helen Li \footnote{Department of Mathematics and Statistics, University of North Carolina at Charlotte,  xli47@uncc.edu} and Molei Tao\footnote{School of Mathematics, Georgia Institute of Technology, mtao@gatech.edu} }

\maketitle

\begin{abstract}
We introduce a novel approach for decomposing and learning every scale of a given multiscale objective function in $\mathbb{R}^d$, where $d\ge 1$. This approach leverages a recently demonstrated implicit bias of the optimization method of gradient descent  \cite{NEURIPS2020_Tao}, which enables the automatic generation of data that nearly follow Gibbs distribution with an effective potential at any desired scale. One application of this automated effective potential modeling is to construct reduced-order models.
For instance, a deterministic surrogate Hamiltonian model can be developed to substantially soften the stiffness that bottlenecks the simulation, while maintaining the accuracy of phase portraits at the scale of interest. Similarly, a stochastic surrogate model
can be constructed at a desired scale, such that both its equilibrium and out-of-equilibrium behaviors (characterized by auto-correlation function and mean path) align with those of a damped mechanical system with the original multiscale function being its potential.
The robustness and efficiency of our proposed approach in multi-dimensional scenarios have been demonstrated through a series of numerical experiments.
A by-product of our development is a method for anisotropic noise estimation and calibration.
More precisely, Langevin model of stochastic mechanical systems may not have isotropic noise in practice, and we provide a systematic algorithm to quantify its covariance matrix without directly measuring the noise. In this case, the system may not admit closed form expression of its invariant distribution either, but with this tool, we can design friction matrix appropriately to calibrate the system so that its invariant distribution has a closed form expression of Gibbs.
\end{abstract}

\section{Introduction}\label{sec:intro}
Physical and data sciences encompass numerous problems that involve multiple spatial and temporal scales. For example, these problems may involve the evolution of a meta-particle within a multiscale energy landscape.

Consequently, a desirable approach to address these challenges is to reduce the complexity of the model in order to improve interpretability, computational efficiency, and reliability. This can be achieved by quantifying the model's features at different scales based on practical interests and separating macroscopic scales from microscopic scales. For instance, it is always desirable to have a surrogate model that captures only the macroscopic behaviors of the full, complex model. This paper proposes an automated and systematic method for decomposing and learning effective models at an arbitrary scale specified by the user, for a given multiscale objective or energy function in multiple dimensions.

Powerful methods for scale decomposition have been developed over several decades. Many of these approaches aim to decompose the scales of a time-series signal or eliminate high-frequency component and/or noise from the time-series \cite{SigProcess_Oppenheim1989,Strang1996,Feng2017,DECHEVEIGNE2019280,Haar1910,Strang1989,Daubechies1992,Cohen1995,Ivan_TimeDecomp2007,Kalman1993,Chou1994,laskar1993frequency,dylewsky2019dynamic}. These are in some sense data-based, model-free approaches and thus are different from ours in scope.

There is also a vast body of literature dedicated to addressing multiscale problems and obtaining effective solutions or surrogate systems at specific scales of interest. Several reviews and books on multiscale modeling and simulations cover this area extensively \cite{attinger2004multiscale,gorban2006model,blanc2006variante,pavliotis2008multiscale,weinan2011principles}. For example, by employing appropriate numerical homogenization and discretization techniques, it is possible to directly obtain effective approximations of solutions at a chosen scale, even without explicitly obtaining the effective systems that produce these solutions. Various methods exist, including but not limited to the equation-free methods \cite{MR2041455, kevrekidis2004equation, kevrekidis2009equation}, the Heterogeneous Multiscale Methods \cite{HMMoriginal, engquist2005heterogeneous,weinan2007heterogeneous,abdulle2012heterogeneous}, FLAVORs \cite{tao2010nonintrusive, FLAVORPDE11}, multigrid methods \cite{bramble1991analysis,bramble2000analysis,lions2001resolution}, generalized finite element methods \cite{babuvska1990p,strouboulis2000design,houston2003sobolev,allaire2005multiscale,gallistl2021mixed}, multiscale Galerkin methods \cite{hou1997multiscale,hughes2006multiscale,wang2008discontinuous,efendiev2009multiscale,cockburn2009analysis}, operator-adapted wavelets \cite{owhadi2017multigrid,owhadi2019operator}, asymptotic-preserving methods \cite{jin2022asymptotic,jin2023asymptotic}, etc. Meanwhile, modal decomposition techniques, such as the reduced basis methods with proper orthogonal decomposition \cite{chorin2007problem,quarteroni2014reduced,hesthaven2016certified} and the dynamic-mode decomposition \cite{schmid2010dynamic,kutz2016dynamic}, provide effective scale decomposition of systems. These breakthroughs mainly focus on approximating the solution, as opposed to constructing approximations of the \emph{system} that produces the solution, although we note there are also great work that construct effective surrogate models simultaneously with or after the solution approximation 
\cite{hou1997multiscale,allaire2005multiscale,efendiev2009multiscale,quarteroni2014reduced,hesthaven2016certified,kutz2016dynamic}.

Instead of directly procuring effective solutions, numerous methods aiming at deriving surrogate models of multiscale problems 
have also been proposed. Here, the surrogate/effective models mainly refer to surrogate differential equations which are simpler and easier to simulate/analyze. They can, for example, be attained through averaging/homogenization \cite{spagnolo1968sulla,de1975sulla,papanicolaou1979boundary,kozlov1979averaging,blanc2006variante,bensoussan2011asymptotic,pavliotis2008multiscale,aizenman2015random,abdulle2012heterogeneous,murat2018h}. The estimation of effective models can also be achieved via spectral information of the underlying systems  \cite{Kessler1999,abdulle2021drift,Abdulle2022}; or through filtering approaches and parameterized kernels  \cite{Kalman1993,Chou1994,GaussianFilter,OPFER2006,majda2012filtering,Griebel_RKHS2015,ZHOU2016}, among others.
Moreover, machine-learning approaches have become increasingly popular for constructing effective models from data. Inference learning can be conducted through statistical tools, such as maximum likelihood, methods of moments, and Bayesian inference \cite{Stuart2007,Krumscheid2013,Abdulle2020a}, or through neural network and deep learning techniques \cite{WeinanE2018,han2019uniformly,alber2019integrating,karniadakis2021physics,peng2021multiscale}.

However, many existing methods are based on analytical tools (e.g., asymptotic analysis such as averaging or homogenization) which assume an explicit scale separation.
Machine learning-based methods may be less constrained by scale separation, but they typically require a good amount of training data {\it a priori}, making it challenging to simultaneously and flexibly learn different scales of interest.

This work proposes a novel approach that complements the existing literature, by providing a method that can decompose and learn any scale of a given multiscale objective function, in multi-dimensional Euclidean spaces. Moreover, in cases where there is only moderate or insignificant scale gaps, it can still construct some effective potential that leads to good approximations (see e.g., Sec.\ref{sec:noScaleSeparation} for more precise statements).

  The main idea of our method is leaned on an interesting fact that smaller scale components of a function can be traded for stochasticity, if one introduces an artificial step of optimizing it using gradient descent, and the scale at which this trading starts can be controlled by the learning rate of gradient descent. This idea is inspired by a recent advance in machine learning, where the problem of whether/how local minima of the training objective function can be escaped is of vital importance. Kong and Tao \cite{NEURIPS2020_Tao} demonstrated how gradient descent with \emph{large} learning rate can quantitatively lead to such escapes, providing an alternative to the common escape mechanics based on noises from stochastic gradients. More precisely, the authors proved that when an objective function exhibits multiscale behaviors and is optimized by gradient descent with a large learning rate (a.k.a. time-step), the deterministic optimization dynamics acts like Langevin dynamics with potential that  only describes the macroscopic part of the objective function, while its microscopic part effectively gets turned into noise via chaotic dynamics. As a result, gradient descent does not converge to a local minimizer of the objective, but instead to a statistical distribution characterized by the macroscopic component of the objective function.

Building upon this insight that a large time-step converts under-resolved scales into noise, we design a self-learning approach to learn components of a function at different scales. Firstly, we artificially introduce a damped deterministic mechanical system, specifically noiseless Langevin dynamics, with the given function as its potential. Then, we choose an appropriate time step size, corresponding to the cutoff scale above which we'd like to approximate this function, and simulate the damped mechanical system using this (large) step size. Next, we simulate the system numerically for a long time, and collect position-momentum values at different time points into a set. Points in this set will approximate a statistical distribution governed by a coarse-grained effective energy function, which however may not be Gibbs distribution yet, due to the an-isotropicity of effective noise. We will thus provide a way to estimate the covariance of the effective noise, which is nontrivial because we don't have access to the effective noise as it is only part of the underlying theory.  Consequently, we will calibrate the damped mechanical system by choosing its friction matrix according to the noise covariance, and simulate the system for a long time again. This time, the collected values will approximately follow Gibbs distribution. Finally, we fit from the data the corresponding potential, which will give the effective approximation of the given function above (including) the designated scale. The resulting function can subsequently be used to construct surrogate models based on practical interests.

\subsection{Problem setup and overview}
\begin{figure}[htp!]
    \centering
    \includegraphics[width = 0.5\textwidth ]{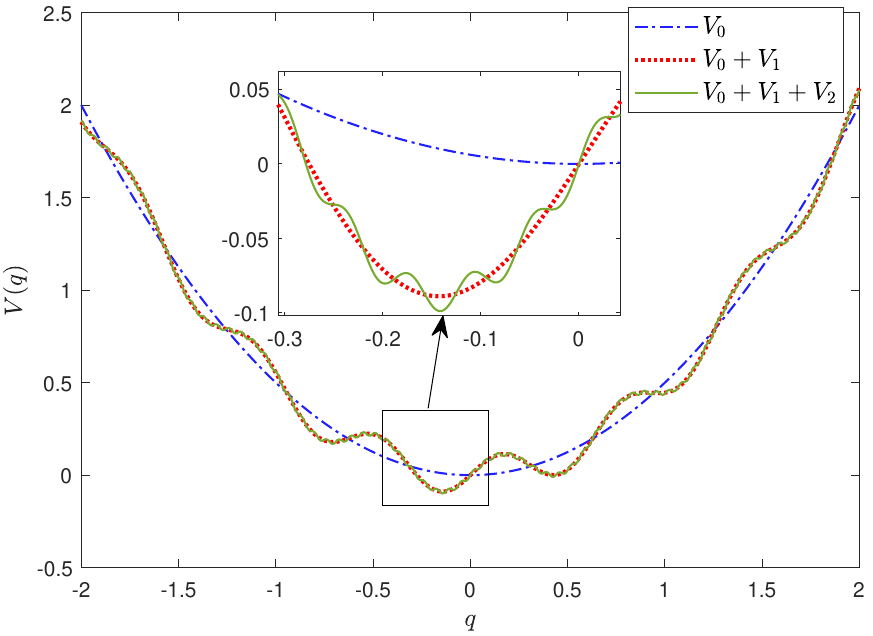}
    \caption{The illustration of a one-dimensional potential function $V=V_0(q)+V_1(q)+V_2(q)$ with three scales: $V_0=\frac{q^2}{2}\sim O(1)$, $V_1=0.1\times \sin(\frac{q}{0.1})\sim O(0.1)$ and $V_2=0.01\times \sin(\frac{q}{0.01})\sim O(0.01)$.
    }
    \label{fig:V_multiDecomp}
\end{figure}

Consider a function $V(\mbf{q})$, which could either be the objective function of an optimization problem, as often encountered in machine learning, or in a conservation law setup, the energy function used for various physical and chemical simulations. Due to the complex nature of these application problems, $V(\mbf{q})$ usually involves multiple scales. By decomposing these scales of $V(\mbf{q})$, we can create a hierarchy of problems that are usually simpler to simulate, analyze, and interpret.
However, it often occurs in practise that we are only given the expression of $V$ without knowing the details of an explicit decomposition. Mathematically, we can thus assume that $V(\mbf{q}): \mathbb{R}^d\rightarrow \mathbb{R}$ with $d\ge 1$ {\it implicitly} admits a multi-scale decomposition
\begin{equation}\label{multiscale_V}
V(\mbf{q})= V_{0}(\mbf{q})+\sum_{j=1}^{K}V_{j,\,\varepsilon_j}(\mbf{q}),
\end{equation}
where $V_0(\cdot)$ denotes the macroscopic component, and other meso- and micro-scale components $\{V_{j,\varepsilon_j}(\cdot)\}_{j=1}^K$ satisfy
$V_{j,\varepsilon_j}=O(\varepsilon_j)$ with scaling parameters $1\gg \varepsilon_1\gg \dots \gg \varepsilon_K>0$. An illustration example of $V(q)$ with three distinct scales is given in Figure~\ref{fig:V_multiDecomp}.

The main goals and contributions of this work are:
\vspace{-0.2 mm}
\begin{itemize}
    \item Our aim is to learn each component of $V$ at different scales, specifically $V_0$ and the series $\{V_{j,\varepsilon_j}\}_j$. We make the assumption that we do not possess explicit access to each individual scale, meaning the values of $\varepsilon_j$' are unknown {\it a priori}, let alone any $V_{j,\varepsilon_j}$.
\item
Once an effective potential function $U_k(\mbf{q})$ that approximates $V$ up to scale $\varepsilon_k$ is learnt (e.g., if $\varepsilon_k \gg \varepsilon_{k+1}$, then $U_k:=V_0+\sum_{j=1}^k V_{j,\varepsilon_j}$, but note $V_{j,\varepsilon_j}$'s are not a priori known), various surrogate models can be constructed, and two will be focused on. In a deterministic case, a surrogate model of Hamiltonian dynamics can produce approximate phase space dynamics with less stiffness and thus allowing much larger time steps for numerical integration. In a stochastic case, a surrogate model can nearly reproduce the equilibrium distribution, auto-correlation function, and mean path of data produced by kinetic Langevin dynamics using the original $V(\mbf{q})$ at the desired scale.
  \item  In order to enable the above contributions, we have also developed an approach to estimate and calibrate kinetic Langevin systems under \emph{an}isotropic stochastic forcings, which pose a challenge as the system's invariant distribution no longer admits an analytical expression, unlike in the isotropic case where the invariant distribution is Gibbs. This is based on estimating the covariance of noise, however without directly measuring it. Instead, we use equilibrium properties only, and by intelligently adjusting the system's mass matrix, a robust and accuracy estimation can be obtained. We can then adjust the system's friction coefficient so that the invariant distribution is put back to the analytically available Gibbs.
\end{itemize}
The organization of this paper is as follows: In Section~\ref{sec:backgound}, we review the theoretical foundations of proposed approach. In Section~\ref{sec:algorithm}, we present the algorithm and explore the analytical attributes of the proposed method. The 1D case is easier as noise will always be isotropic, and it is first discussed. Then we detail the workings of the algorithm in the multivariate  scenario, after describing our additional tool for kinetic Langevin estimation and calibration. Then, Section~\ref{sec:numerics} delivers numerical simulations of various examples in both 1D and higher dimension. These results substantiate the robustness and efficiency of our proposed approach. Section~\ref{sec:conclusion} summarizes this study and suggests possible future directions.

\section{Motivation and theoretical background}\label{sec:backgound}
In this section, we will briefly review some theoretical results in \cite{NEURIPS2020_Tao} that motivate, and will be used by this work. We will also review some background knowledge about kinetic Langevin  which will be used later.

\subsection{Large step size effectively turns under-resolved scale into noise}
\label{sec:reviewKongTao}
\paragraph{Informal version} For quantitative understanding of deep learning, Kong and Tao \cite{NEURIPS2020_Tao} quantified how the optimization algorithm of gradient descent escapes local minimizers of a function. They considered the optimization of a multiscale objective function by gradient descent with large learning rate, i.e. iterations given by
\begin{equation}
    x_{k+1}=x_k-\eta \nabla f(x_k),
    \label{eq:gradientDescent}
\end{equation}
where $f:\mathbb{R}^d\rightarrow \mathbb{R}$ was assumed to be admitting a macro-micro decomposition given by
\[
f(x): = f_0(x)+f_{1,\varepsilon} (x)
\]
where $f_0$ is the macroscopic, $\mathcal{O}(1)$ component of the potential, and $f_{1,\varepsilon}(x)=\epsilon f_1(x/\epsilon)$ is the microscopic component as it squeezes an $\mathcal{O}(1)$ function $f_1$ in both the x and y directions. $\eta$ was referred to as the learning rate. Large learning rate is in the sense that $\eta\gg 1/L$, where $L$ is the Lipschitz constant of $\nabla f$, because classical analysis of gradient descent assumed $\eta < 1/L$, in which case the iterates provably converge to a (nearby) local minimizer of $f$.

They showed that in this case, the iteration $x_{k+1}=x_k-\eta \nabla f_0(x_k) - \eta \nabla f_1(x_k/\epsilon)$, where albeit being deterministic, exhibits behaviors similar to those of a stochastic iteration
\begin{equation}
    x_{k+1}=x_k-\eta \nabla f(x_k) - \eta \zeta_k,
    \label{eq:GDauxiliaryDynamics}
\end{equation}
where $\zeta_k$'s are random variables i.i.d. to $\zeta$ which will be defined. The condition $\eta \gg 1/L$ actually leads to $\eta \gg \epsilon$, which means step size $\eta$ is too large to resolve the details of $f_{1,\epsilon}$. As a consequence, the effective contribution from the microscopic, under-resolved gradients $\nabla f_1(x_k/\epsilon)$ over many iterations can be approximated by contributions from noise $\zeta_k$. Therefore, gradient descent \eqref{eq:gradientDescent} with \textbf{large} learning rate $\eta$ is \textbf{not} a time discretization of gradient flow ODE
\[
    \dot{x}=-\nabla f(x) = -\nabla f_0(x) - \nabla f_1(x/\epsilon)
\]
but can be intuitively understood as a stepsize $\eta$ time discretization of \textbf{SDE}
\begin{equation}
    dx= -\nabla f_0(x) dt + \sqrt{\eta} \Sigma dW_t
    \label{eq:GDauxiliarySDE}
\end{equation}
for some constant covariance matrix $\Sigma$ given by $\Sigma \Sigma^T = \mathbb{E} [\zeta \zeta^T]$.

\paragraph{Selected rigorous details} They considered $\epsilon\ll 1$, $f_0,\,f_{1,\epsilon}\in \mathcal{C}^2(\mathbb{R}^d)$ and microscopic potentials $f_{1,\epsilon}$ satisfying the following two conditions:
\begin{itemize}
\item {\bf $\mathcal{O}(1)$ gradient condition: } there exists a bounded and nonconstant random variable (r.v.) $\zeta\in \mathbb{R}^d$ with $\mathbb{E}(\zeta)=0$, such that: $\forall \epsilon>0$ and $\forall x\in \mathbb{R}^d$ there exists a positive measured set $\Gamma_{x,\epsilon}\subset B(x, \delta(\epsilon)$ with $\lim_{\epsilon \downarrow 0}\delta (\epsilon)=0$, such that a uniformly distributed r.v. on $\Gamma_{x,\epsilon}$, denoted by $u_{\Gamma_{x,\epsilon}}$, satisfies $\nabla f_{1,\epsilon}(x+u_{\Gamma_{x,\epsilon}})\xrightarrow{w}-\zeta$ as $\epsilon\rightarrow 0$, uniformly with respect to $x$.
\item {\bf $\mathcal{O}(1/\epsilon)$ Hessian condition: }$\epsilon \nabla^2 f_{1,\epsilon}$ is uniformly bounded as
$\epsilon \rightarrow 0$. Also $\exists m>0$, such that for any bounded and nonempty rectangle $\Gamma\subset \mathbb{R}^d$, the expectation  $\mathbb{E}\left[\ln \|\epsilon \nabla^2 f_{1,\epsilon} (u_{\Gamma})\|_2\right]\rightarrow m$ as $\epsilon\rightarrow 0$, where $u_{\Gamma}$ is a uniform r.v. on $\Gamma$.
\end{itemize}
While these conditions are nontrivial, Kong and Tao noted that they include but are strictly weaker than periodicity (e.g., $f_{1,\epsilon}(x)=\epsilon  \sin(x/\epsilon)$, for which $\Sigma=1/\sqrt{2}$) or quasiperiodicity.

They then considered the deterministic dynamical system induced by iteration map
\begin{equation}\label{GD}
\varphi: x\mapsto x-\eta \nabla f(x) = x-\eta \nabla_0 f(x)-\eta \nabla f_{1,\epsilon}(x),
\end{equation}
and proved the following results (selected):
\begin{enumerate}
\item Fix $\eta$ and let $\epsilon\to 0$. If $\varphi$ has a family of nondegenerate invariant distributions for $\{\epsilon_i\}_{i=1}^{\infty}\to 0$ which converges in the weak sense, then the weak limit is an
invariant distribution of $\hat{\varphi}$, where $\hat{\varphi}$ defines a stochastic map
\begin{equation}\label{SGD}
\hat{\varphi}: x\mapsto x-\eta \nabla f_0(x)+\eta \zeta.
\end{equation}
\item In a special case where the macroscopic potential $f_0$ is strongly convex and $L$-smooth, there exists some constant $C>0$ independent of $\epsilon$, such that the mapping $\hat{\varphi}$ has a unique invariant distribution for any fixed $\eta\le C$ and the iteration of $\hat{\varphi}$ converges exponentially fast to this distribution. Moreover, if the covariance matrix of $\zeta$ is isotropic, i.e., $\sigma^2\mathrm{I}_d$, then the rescaled Gibbs distribution $\frac{1}{Z}\exp\left(-\frac{2f_0(x)}{\eta \sigma^2}\right)\, \mathrm{d}x$ is an $O(\eta^2)$ approximation of that of $\hat{\varphi}$.
\end{enumerate}
However, it is worth noting that \cite{NEURIPS2020_Tao} only considered two scales in the objective function $f(x)$, whereas this work heuristically extends the application to arbitrary $K$ scales. In addition, note $\zeta_k$ in \eqref{eq:GDauxiliaryDynamics} is actually bounded, unlike `$dW_t$' in \eqref{eq:GDauxiliarySDE}, but \eqref{eq:GDauxiliarySDE} is still a reasonable interpretation in a central limit theorem sense (see \cite{NEURIPS2020_Tao} for details). Additionally, the stochastic behavior of the deterministic mapping with a large LR was only proved for gradient descent dynamics, and if one considers a damped deterministic mechanical system instead, with $\nabla f$ being its forcing, and discretize its time using a stepsize that under-resolves part of $f$, the quantitative stochasticity of its dynamics was only a conjecture.
While this work relies on this conjecture, our empirical results appear to be consistent with it, namely deterministic simulations of damped mechanical systems with large step size exhibit behaviors resemblant to (stochastic) kinetic Langevin dynamics.

\subsection{Review of kinetic Langevin dynamics}\label{subsec:Langevin}
Kinetic Langevin equation models how  particles evolve subject to a combination of conservative forcing, damping, and thermal noise. Let $(\mbf{q},\,\mbf{p})$ denote the position and momentum of particles, and $M$ denote the corresponding mass matrix, then it gives the evolution of $(\mbf{q},\,\mbf{p})$ by
\begin{equation}
\begin{cases}
d \mbf{q}= M^{-1}\mbf{p} dt ,\\
d \mbf{p} = (-\nabla V(\mbf{q})-\Gamma M^{-1}\mbf{p}) dt + \Sigma d\mbf{W}_t,
\end{cases}
\label{eq:kineticLangevin}
\end{equation}
where $\nabla V$ denotes the potential force, $\mbf{W}_t$ denotes the standard multivariate Wiener process, $M$ is a positive definite mass matrix, $\Gamma$ is the friction coefficient matrix for the dissipative force, and $\Sigma$ denotes the diffusion coefficient matrix for the thermal random forces.

It is known (e.g., \cite{Pavliotis}) that if the diffusion coefficient matrix $\Sigma$ and the friction coefficient matrix $\Gamma$ satisfy the following fluctuation–dissipation relation
\begin{equation}\label{FDR}
\Sigma \Sigma^T =2\Gamma,
\end{equation}
then the equilibrium distribution is a Gibbs distribution
\begin{equation}\label{equib_dist}
d\pi(\mbf{q},\,\mbf{p})
=\frac{1}{Z}\exp{-(\left(V(\mbf{q})+\mbf{p}^TM^{-1}\mbf{p}/2)\right)} d\mbf{q}d\mbf{p}.
\end{equation}
The most commonly discussed case is when the thermal noises are isotropic, that is $\Sigma \Sigma^T =\sigma^2\mathrm{I}_d $ for some scalar $\sigma$. In this case, the equilibrium distribution is always Gibbs as long as $\Gamma=\frac{\sigma^2}{2}\mathrm{I}_d =\gamma \mathrm{I}_d$.

However, if the noises are anistropic, the invariant distribution of \eqref{eq:kineticLangevin} still exists under mild conditions, but it may not admit an analytical expression.

\section{Main results}\label{sec:algorithm}
In this section, we propose a self-learning algorithm that automatically identifies macroscopic scales of a function $V$ (eq.\ref{multiscale_V}) up to a desired scale
$\varepsilon$. Here, $\varepsilon$ corresponds to a cutoff level $k$ with  $1\le k\le K$ satisfying $\varepsilon_k \gtrsim \varepsilon \gg \varepsilon_{k+1}$, meaning we seek $U_k=V_0+\sum\limits_{j=1}^{k}V_{j,\,\varepsilon_j}$.
Of course, if the goal is to separate all the scales instead, we can just repeat this procedure for different $\varepsilon$ values.

\subsection{The proposed self-learning methodology for multi-dimensions}
Our \underline{first step} of identifying components of $V$ above scale $\varepsilon$, solely based on $V$\footnote{More precisely, only 1st-order oracle, no additional data or information.}, is to extend the gradient descent iterations analyzed in \cite{NEURIPS2020_Tao} (see Sec.\ref{sec:reviewKongTao} for a quick summary) by including momentum.

More precisely, we consider a large step size simulation of the following damped mechanical system
\begin{equation}\label{Newton_sys}
\begin{cases}
\dot{\mbf{q}}= M^{-1}\mbf{p},\\
\dot{\mbf{p}}=-\nabla V(\mbf{q})-\Gamma M^{-1}\mbf{p},
\end{cases}
\end{equation}
where the step size is chosen to be at the order of $\varepsilon$, i.e. $\delta\sim \varepsilon$, so that it underresolves undesired smaller scales. Here $M$ and $\Gamma$ are mass and friction coefficient matrices that we can choose, but $V(\mbf{q})$ is given and we have to use it in its entirety as its scale decomposition is not yet known (this actually poses a significant challenge in high dimension in general, which will be described and addressed in Sec.\ref{sec:findSigma}).

The choice of numerical integrator for \eqref{Newton_sys} for the purpose of turning unresolved small scales effectively into `noise' needs not to be unique, but here we simply apply a dissipative generalization of St\"ormer--Verlet, where friction is handled by an exponential integrator. The evolution of discrete solution $(\mbf{q}_n,\, \mbf{p}_n)$ at time grid $n$ thus becomes
\begin{equation}\label{verlet_scheme}
    \begin{aligned}
    \mbf{q}_{n+\frac{1}{2}}=&\mbf{q}_n+M^{-1}\frac{\delta}{2} \mbf{p}_{n},\\
    {\mbf{p}}_{n+1}=& e^{-\Gamma M^{-1} \delta}\,\mbf{p}_{n}-\delta \nabla V(\mbf{q}_{n+\frac{1}{2}}),\\
    \mbf{q}_{n+1}=&\mbf{q}_{n+\frac{1}{2}}+M^{-1}\frac{\delta}{2} {\mbf{p}}_{n+1},
    \end{aligned}
\end{equation}
and with given initial condition $(\mbf{q}_0,\, \mbf{p}_0)$ from a suitable region. We conjecture that the gradient of under-resolved microscopic components of $V$ can effectively be approximated by noises as $\varepsilon_{k+1}\to 0$, similar to Sec.\ref{sec:reviewKongTao}. This conjecture means that invariant distribution of the deterministic dynamics  can be approximated by that of the stochastic system
\begin{equation}
\begin{cases}
d \mbf{q}= M^{-1}\mbf{p} dt ,\\
d \mbf{p} = (-\nabla U_{k}(\mbf{q})-\Gamma M^{-1}\mbf{p}) dt + \Sigma dW_t,
\end{cases}
\label{eq:kineticLangevin_v2}
\end{equation}
with the effective potential $U_{k}=V_0+\sum_{j=1}^{k}V_{j,\varepsilon_j}$, $\Gamma=\frac{1}{2}\Sigma\Sigma^T$, and the diffusion coefficient matrix given by
\begin{equation}\label{mico_variance}
\begin{aligned}
\Sigma\Sigma^T :
&=
 \delta \mathbb{E}\left[\nabla \left(\sum_{j=k+1}^{K} V_{j,\varepsilon_j}(\mbf{u})\right)\otimes \nabla \left( \sum_{j=k+1}^{K} V_{j,\varepsilon_j}(\mbf{u}) \right)\right],
 \end{aligned}
\end{equation}
where the expectation is with respect to an auxiliary random variable $\bm u$ defined in the gradient conditions in Section~\ref{sec:backgound}, i.e. a uniformly distributed r.v. on any bounded rectangle in $\mathbb{R}^d$ with size independent of $\varepsilon_{k+1}, \cdots, \varepsilon_{K}$, i.e. $\mathcal{O}(1)$. The sense of approximation is that the invariant distribution of \eqref{eq:kineticLangevin_v2} differs from that of \eqref{verlet_scheme} by at most $\mathcal{O}(\delta)$ (in weak* topology) as scale separation goes to infinity (i.e. $\delta/\varepsilon_{k+1} \to \infty$).

  In addition, if the system is mixing (unfortunately, the precise necessary and sufficient condition for so is still a major open problem, but it was conjectured to be the case for large learning rate gradient descent;
  see Section~\ref{sec:backgound}
  for semi-quantitative discussions), then \eqref{verlet_scheme} will converge to its invariant/equilibrium distribution, denoted by $\pi_\delta$.
  This in the sense that 1) an ensemble of trajectories with random initial condition following any smooth density will converge to the invariant distribution as $n\to\infty$, and 2) except for measure zero initial conditions, any single trajectory will have an ergodic limit with respect to the invariant distribution $\pi_\delta$ , meaning for any smooth test function $\phi$, we have
  \[
    \lim_{N\to\infty} \frac{1}{N}\sum_{n=1}^N \phi(\bm q_n, \bm p_n) = \mathbb{E}_{\bm q, \bm p \sim \pi_\delta} \phi(\bm q, \bm p).
  \]
This means, given data $\{\mbf{q}_n,\mbf{p}_n\}$ generated from $V$ via \eqref{verlet_scheme}, the empirical distribution of $\{\mbf{q}_n,\mbf{p}_n\}_{n=1,\cdots, N}$ converges to $\pi_\delta$, which is an $\mathcal{O}(\delta)$ approximation of $\pi$, whose density is $Z^{-1} \exp(-U_k(\bm q)-\bm p^T M^{-1} \bm p/2)$ for \eqref{eq:kineticLangevin_v2}.

Therefore, for a given $\Gamma$, \underline{Step Two} of our approach is to run the deterministic system \eqref{verlet_scheme} for sufficiently long on time interval $[0, T]$ until $(\mbf{q}_n,\, \mbf{p}_n)$ reach their statistical equilibrium, which can be tested from the normality test \cite{thode2002testing} of the distribution of $\{\mbf{p}_n\}$. Note our introduction of momentum is rather important, despite the fact that the whole dynamics is just artificially introduced, because the convergence of $\mbf{p}$ to simple Gaussian enables accurate detection of near convergence.

This Step Two generates a training data set solely out of $V$, which will be used later for the learning of an effective macroscopic function, hence the name `self-learning'.

\begin{remark}\label{rmk:init_cond}
Although in principle it is very flexible to choose the initial condition $(\mbf{q}_0,\, \mbf{p}_0)$ because of the stochasticity induced by microsuitable large step size, it becomes subtle in practice as good initial conditions can greatly improve the explorcity and speed-up convergence of the algorithm. We will provide more details for specific examples in Section~\ref{sec:numerics}.
\end{remark}

\begin{remark}\label{rmk_unkown_eps}
We explain the strategy here when the scale separation among $\{\varepsilon_j\}_{j=1}^K$ exists with range of scales knowing, but the value of each $\varepsilon_j$ keeps unknown.

The idea is to firstly use an isotropic friction matrix $\Gamma=\gamma I$ and run the damped mechanical system  \eqref{verlet_scheme} using $\delta$ from an array of time-stepping sizes, which are arranged in a decreasing order and are compatible with the range of scales of $V(\mbf{q})$.
For each selected $\delta$,
we simulate long enough until  numerically reach the equilibrium distribution of $\pi_{{\delta}}(\mbf{q},\,\mbf{p})$. If the distribution of $\mbf{p}$ becomes a general multi-dimensional normal distribution via a normality test, then it suggests we find a scale $\varepsilon_k$ which is proportional to $\delta$. Otherwise, if the distribution of $\mbf{p}$ is far away from a normal distribution, then this suggests the system is still under mixed scales and we will further decease $\delta$ and regenerate data. After we find the proper $\delta$ at a targeting scale, we next tune the friction matrix $\Gamma$  to achieve the Gibbs distribution \eqref{Gibbs_dist}.

The successful results of a benchmark test with unknown scales can be found in Figure~\ref{Fig:cosHighF_Potlearnt}. Also, notice that the simulations on different $\delta$ are parallelizable, so the total efficiency of distinguish scales will be ensured.
\end{remark}

\subsection{Microscopic covariance estimation and the creation of fluctation-dissipation balance}
\label{sec:findSigma}

Recall that the friction coefficient $\Gamma$ is a parameter we can choose, while $\Sigma$ is something we cannot control because it characters the effective noise that originates from the microscopic gradient. If we can align $\Gamma$ with $\Sigma$ so that they commute, or more precisely so that $\Sigma \Sigma^T = 2\Gamma/\hat{\beta}$ for some constant scalar $\hat{\beta}$ (known as the inverse temperature), then the data collected in the above Step 2 will approximately follow a density $\propto \exp{\left(-\hat{\beta}\left(U_k(\mbf{q})+\mbf{p}^TM^{-1}\mbf{p}/2\right)\right)}$, which means we can recover $U_k$ via regression of the $\mbf{q}$ marginal of the data.

However, we do not know $\Sigma$ a priori because we only have $V$, but not the microscopic potential $U_k$, upon which $\Sigma$ is based. This thus poses a challenge.

When $V$ is one-dimensional, this challenge is relatively easy to solve, as aligning $\Gamma$ with $\Sigma$ is not an issue and we just need the inverse temperature, which can be estimated from $\mathbf{p}$'s variance. When $V$ is a multi-dimensional function in general, choosing $\Gamma$ so that the system admits Gibbs as its invariant distribution is rather nontrivial, but we will provide a rigorous and systematic method. Both cases are now detailed:

\paragraph{One-dimensional strategy}
When $d=1$, the friction matrix $\Gamma$ becomes just a scalar $\gamma$. Hence, the equilibrium distribution of $(q,p)$ is always a Gibbs distribution. 
Consequently, $M$ and $\gamma$ can be selected easily in \eqref{Newton_sys}. Throughout the paper, we select the following values for the 1D simulation $0<\gamma \le 0.1\text{ and }  M = 1.$ Once the numerical solution $\{({q}_n, {p}_n)\}_n$ is nearly in the statistical equilibrium, the distribution of ${p}_n$ will be close to a Gaussian distribution. Following \cite{Pavliotis14}, we estimate the equivalent thermal noise effect $\hat{\beta}$ induced by the gradient of micro-scale components
\begin{equation}\label{equiv_temp}
\hat{\beta}=\frac{1}{\mathsf{var}(p)}.
\end{equation}
In 1D, $\pi_{\delta} (q,\,p)$ is approximated by
\begin{equation}\label{Gibbs_dist}
\pi_{\delta} (q,\,p)\sim \exp\left(-\hat{\beta}\big(U_{k}(q)+\frac{p^2}{2}\big)\right)dq\, dp.
\end{equation}
In 1D simulation with suitable value of $\gamma$, we can get $\hat{\beta}$ in the range around $\hat{\beta}\in (0.2,\, 1.25)$, which will lead to better learning results for $U_k$.

\paragraph{Multi-dimensional strategy}
As mentioned above, due to the anisotropy of microscopic scales in $V(\mbf{q})$ the dimension $d\ge 2$, invariant distribution of kinetic Langevin \eqref{eq:kineticLangevin_v2} in general may not be interpretable. Therefore, we construct a special case by estimating 
$\Sigma$ (originated from the microscopic scales \eqref{mico_variance}) and choose 
$\Gamma=\frac{1}{2}\Sigma \Sigma^T$ to ensure the fluctuation-dissipation relation \eqref{FDR}. In this case, the temperature is always $\hat{\beta}=1$. We explain the details by a two-stage scheme.
\noindent
\paragraph{\bf Stage 1: estimate $\Sigma\Sigma^T$.}
The key idea of our innovation is to exploit the effects of different choices of mass matrix $M$ in
the following stochastic Langevin dynamics:
\begin{equation}\label{Langevin_sys2}
\begin{cases}
{d\mbf{q}} =& M^{-1}\mbf{p}dt,\\
{d\mbf{p}} = & \left(-\gamma M^{-1} \mbf{p}-\nabla V(\mbf{q})\right)dt + \Sigma d\mbf{W}_t.
\end{cases}
\end{equation}
where $0<\gamma$ is some scalar friction constant. Note $\nabla V dt + \Sigma d\mbf{W}_t$ come together from the under-resolved small scales, and thus again $\Sigma$ cannot be chosen or determined a priori. The positive-definite mass matrix $M$ however can be arbitrarily chosen, because \eqref{Langevin_sys2} is just some auxiliary dynamics designed for obtaining an algorithm. We will design $M$ values so that $\Sigma$ can be estimated from `solutions' of \eqref{Langevin_sys2}, which are actually collected as 
the deterministic iterations with step size $\delta$
\[
\begin{aligned}
    \mbf{q}_{n+\frac{1}{2}}=&\mbf{q}_n+M^{-1}\frac{\delta}{2} \mbf{p}_{n},\\
    {\mbf{p}}_{n+1}=& e^{-\gamma M^{-1} \delta}\,\mbf{p}_{n}-\delta \nabla V(\mbf{q}_{n+\frac{1}{2}}),\\
    \mbf{q}_{n+1}=&\mbf{q}_{n+\frac{1}{2}}+M^{-1}\frac{\delta}{2} {\mbf{p}}_{n+1}.
    \end{aligned}
    \]
Meanwhile, the relation between $\Sigma\Sigma^T$ and the solutions is established in the following proposition.

\begin{theorem}\label{thm_Sigma}
Consider a stochastic Langevin dynamics \eqref{Langevin_sys2} with a positive scalar friction constant $\gamma$ and a constant diffusion coefficient matrix $\Sigma$, then $\Sigma$ satisfies
\begin{equation}\label{Sigma_est}
\mathsf{Tr}\left((\Sigma\Sigma^T )M^{-1}\right)=\gamma \mathbb{E}\left[ \mbf{p}^{T}M^{-2}\mbf{p}\right],
\end{equation}
where the expectation is taken with respect to the invariant distribution of $(\mbf{q},\,\mbf{p})$.
\end{theorem}
\begin{proof}
We recall the definition of Hamiltonian $H(\mbf{q},\,\mbf{p})=\frac{1}{2}\mbf{p}^T M^{-1}\mbf{p}+V(\mbf{q})$ and evaluate it along the solutions of \eqref{Langevin_sys2}. It\^o formula shows that the Hamiltonian satisfies the following SDE
\[
\begin{aligned}
d H = & (M^{-1}\mbf{p})^{T} \left[-\gamma M^{-1}\mbf{p}-\nabla V(q)\right]dt +\nabla V(q)^{T} M^{-1}\mbf{p}dt + \mathsf{Tr}\left(( \Sigma\Sigma^T )M^{-1}\right)dt\\
&\quad +\left(M^{-1}\mbf{p}\right)^{T}\sqrt{2}\Sigma d\mbf{W}_t\\
=& -\gamma \mbf{p}^T M^{-2}\mbf{p}dt+
\mathsf{Tr}\left((\Sigma\Sigma^T )M^{-1}\right)dt+\sqrt{2}\mbf{p}^{T}M^{-1}\Sigma d\mbf{W}_t.
\end{aligned}
\]
Integrating it from time $0$ to $t$, we get
\[
\begin{aligned}
H(t)=H(0)+&\int_{0}^{t}-\gamma \mbf{p}^{T}(\tau)M^{-2}\mbf{p}(\tau)d\tau+  \int_0^t \mathsf{Tr}\left((\Sigma\Sigma^T )M^{-1}\right)d\tau\\ &\quad +\sqrt{2}\int_0^{t}\mbf{p}^{T}(\tau)M^{-1}\Sigma d\mbf{W}_\tau.
\end{aligned}
\]
We take expectation of $H$ with respect to the invariant distribution and notice that the last It\^o integral term is a martingale, so
\[
\mathbb{E}(H(t))=:E(t)=E(0)
+ \int_{0}^{t}-\gamma\mathbb{E}\left[ \mbf{p}^{T}(\tau)M^{-2}\mbf{p}(\tau)\right]d\tau+ \int_0^t \mathsf{Tr}\left((\Sigma\Sigma^T )M^{-1}\right)d\tau.
\]
Hence we get the differential equation of $\mathbb{E}(H(t))$, which is
\[
\begin{aligned}
\dot{E}(t)=&-\gamma \mathbb{E}\left[ \mbf{p}^{T}(t)M^{-2}\mbf{p}(t)\right]+
\mathsf{Tr}\left((\Sigma\Sigma^T )M^{-1}\right).
\end{aligned}
\]
When the system reaches its statistical equilibrium, we have $\dot{E}(t)=0$, so
\[
\begin{aligned}
0=\dot{E}(t)=&-\gamma \mathbb{E}\left[ \mbf{p}^{T}(t)M^{-2}\mbf{p}(t)\right]+
 \mathsf{Tr}\left((\Sigma\Sigma^T )M^{-1}\right).
\end{aligned}
\]
Consequently, we can set up a relation between $M$ and $\Sigma$ at the statistical equilibrium of system, that is
\[
\mathsf{Tr}\left((\Sigma\Sigma^T )M^{-1}\right)=\gamma \mathbb{E}\left[ \mbf{p}^{T}M^{-2}\mbf{p}\right],
\]
which proved the proposition.
\end{proof}
\noindent
In order to estimate $Z:=\big(\Sigma \Sigma^T\big)=\left(z_{ij}\right)_{d\times d}$, we then can design the entries of inverse of mass matrix $A:=M^{-1}=\left(a_{ij}\right)$ accordingly and apply the following steps:
\begin{itemize}
\item Step 1: choose a set of $d$ different diagonal matrices  $A^{(k)}=\mathsf{diag}(a_{ii}^{(k)})$ and run the Langevin dynamics \eqref{Langevin_sys2} to reach its statistical equilibrium. Thus, we get a system of $d$ linear equations to solve all diagonal entries $z_{ii}$:
\[
 \sum_{i=1}^{d} a_{ii}^{(k)}z_{ii} =
\gamma \mathbb{E}\left[ {\mbf{p}^{(k)}}^{T}M^{-2}\mbf{p}^{(k)}\right],\quad k=1,\dots, d.
\]
\item Step 2: fix a pair of off-diagonal entries with $\ell\neq r$, set this pair $a_{r\ell}=a_{\ell r}=1/2$, set all diagonal entries $a_{ii}=1$, and set the rest of off-diagonal entries $a_{ij}=0$, then run the Langevin dynamics \eqref{Langevin_sys2} until reaching the equilibrium. We thus have
\[
\left(z_{\ell r}+ \sum_{i=1}^{d} z_{ii}\right) = \gamma \mathbb{E}\left[ \mbf{p}^{T}M^{-2}\mbf{p}\right],
\]
which can solve $z_{\ell r}$. Hence, after repeating this procedure for $\frac{d(d-1)}{2}$ times, we can get all values of the off-diagonal entries of $Z=(\Sigma\Sigma^T)$.
\end{itemize}
Combining step 1 with step 2, overall we need to solve $\frac{d(d+1)}{2}$ linear systems. However, each one of them is independent of the others which means all the simulations can be run in parallel.

\begin{remark}\label{rmk:Sigma_2d}
To demonstrate the scheme in details, we consider a two-dimensional case $d=2$. In step 1, we choose the two inverse of mass matrices to be
\[
A^{(1)}=M^{-1, (1)}=\begin{pmatrix}1 & 0\\ 0 & 2 \end{pmatrix}
\quad \text{ and }\quad A^{(2)}=M^{-1, (2)}=\begin{pmatrix}2 & 0\\ 0 & 1 \end{pmatrix}.
\]
In step 2, we choose the inverse of mass matrices to be
\[
A^{(3)}=M^{-1,(3)}=\begin{pmatrix}1 & 1/2\\ 1/2  & 1\end{pmatrix}.
\]
\end{remark}
\paragraph{\bf Stage 2: estimate $\Gamma$.}
Once we get the estimation of $\Sigma$ from stage 1, we set the mass matrix to be $M=\mathbb{I}_d$ and set the friction coefficient matrix $\Gamma$ to be
\[
2 \Gamma = \Sigma \Sigma^T =Z.
\]

\subsection{Estimate the effective potential}\label{sub:interpolation}
We apply the Verlet scheme \eqref{verlet_scheme} to the damped mechanical system \eqref{Newton_sys} (not kinetic Langevin), with appropriate $\Gamma$, for sufficiently long, and record the trajectory as $\{\mbf{q}_n, \mbf{p}_n\}$, i.e. our data. We propose to learn the effective potential $U_k$
from the invariant distribution \eqref{Gibbs_dist}, which is
\begin{equation}\label{Pi_EffectU}
\log \left(\pi_{\delta}(\mbf{q}) \right)=-\log Z+\hat{\beta} U_k(\mbf{q}),
\end{equation}
where $\pi_{\delta}(\mbf{q})$ is the marginal of invariant distribution, approximated empirically by $\{\mbf{q}_n\}$. For 1D case, $\hat\beta =\frac{1}{\mathsf{var}(p)}$, while for multi-dimensional case $\hat\beta=1$ as $\Gamma$ is set to satisfy $\Gamma = \frac{1}{2}\Sigma \Sigma^T$.

Learning $U_k(\mbf{q})$ from data via \eqref{Pi_EffectU} is a function approximation / interpolation problem. There are of course many ways to interpolate the data, but since function approximation is not the main point of this paper, we will just pick one approach. Specifically, we choose a set of suitable basis functions, for instance the piecewise spline basis,
$\{\varphi_i(\mbf{q})\}_{i}^{m}$ to approximate $U_k(\mbf{q})$ by \begin{equation*}
U_k(\mbf{q})\approx \sum_{i=1}^{m}a_i\varphi_i(\mbf{q}),
\end{equation*}
where we can solve a regression problem to learn the coefficients $\{a_i\}_{i=1}^{m}$.

\subsection{Summary of the algorithm}\label{sub:method_summay}
\noindent
See Algorithm \ref{Algorithm_summary}.
\begin{algorithm}
\caption{Learn effective function $V=U_k+\mathcal{O}(\varepsilon_k)$.
\label{Algorithm_summary}}

\textbf{Input:} Target function $V(\mbf{q})$, a step size $\delta\sim O(\varepsilon_k)$.\\
\textbf{Output:} Effective function $U_k(\mbf{q}) \sim O(\varepsilon_k)$ and effective covariance matrix $\Sigma\Sigma^T$ from unresolved scales.
\begin{algorithmic}[1]
\State {\bf Stage 1 and 2 :} estimate $\Sigma\Sigma^T$ and $\Gamma$: \newline
Set $\Gamma = \gamma I$ with $0<\gamma \le 0.1$, fix the step size $\delta$, run \eqref{verlet_scheme} with $V(\mbf{q})$ for a total of $\frac{d(d+1)}{2}$ independent trajectories using different mass matrices. Solve entries of $\Sigma\Sigma^T$ via \eqref{Sigma_est}. Set $\Gamma = \Sigma\Sigma^T/2$.
\State 
Estimate $\pi_\delta(\mbf{q})$ and $U_k(\mbf{q})$:\newline
Set $M=I$, fix the step size $\delta$, run \eqref{verlet_scheme} with $V(\mbf{q})$ for one trajectory. Estimate the empirical distribution $\pi_\delta (\mbf{q})$ from the data $\{\mbf{q}_n\}$. Learn $U_k(\mbf{q})$ via \eqref{Gibbs_dist}.
\end{algorithmic}
\end{algorithm}

\section{Surrogate models and numerical experiments}\label{sec:numerics}

We now numerically test the efficacy of our approach in two senses. One is about the accuracy of the learned effective function $U_k(\mbf{q})$. The other is about the approximation abilities of derived surrogate models; i.e., if the learned $U_k$ replaces the full $V$ when used in a dynamical setup, how would the resulting dynamics differ. For this latter point, we will illustrate two dynamical setups, one being Hamiltonian dynamics and the other being kinetic Langevin, in both cases $V$ or $U_k$ will be the potential. The versions with $U_k$ are the surrogate models, whose benefits are that they can be simulated using larger setups, easier to analyze, and sometimes more intuitive to interpret as well.

More precisely, the surrogate models are given by
\begin{equation}\label{surr_Hamiltonian}
\text{\bf Hamiltonian: }    \begin{cases}
\dot{\mbf{q}}= M^{-1}\mbf{p},\\
\dot{\mbf{p}}=-\nabla U_k(\mbf{q}),
\end{cases}
\end{equation}
and
\begin{equation}\label{Surr_Langevin}
\text{\bf Kinetic Langevin:} \begin{cases}
\dot{\mbf{q}}= M^{-1}\mbf{p},\\
\dot{\mbf{p}}=-\nabla U_k(\mbf{q})-\Gamma M^{-1} \mbf{p}+\sqrt{2\Gamma} d\mbf{W}_t.
\end{cases}
\end{equation}
For the comparison of the deterministic Hamiltonian system, we will mainly focus on the comparison of the phase portraits $(\mbf{q},\, \mbf{p})$ generated by $V(\mbf{q})$ and $U_k(\mbf{q})$. For the case of kinetic Langevin, we will compare the equilibrium distribution, mean path (ensemble average of trajectories) and normalized auto-correlation functions (Normalized ACF) of $\mbf{q}$.
 Note that this is not a tautology, as we will not only inspect the equilibrium aspect of the surrogate model, which should match that of the full model as long as the proposed trading of small-scale for effective noise is correct, but we will also investigate out-of-equilibrium dynamical aspects by comparing the mean trajectories and normalized ACF.

(Empirical) mean trajectory is defined as
\begin{equation} \label{eq:mean_traj}
\bar{\mbf{q}}(t_n) = \frac{1}{M}\sum_{\ell=1}^M\mbf{q}^{(\ell)}(t_n),
\end{equation}
for which we evolve $M$ independent trajectories from i.i.d. random initial condition, and $\mbf{q}^{(\ell)}(t_n) $ denotes the $\ell$-th one of them.
We also define the normalized ACF as
\begin{align} \label{eq:ACF}
\text{Normalized ACF}(\tau)=\frac{1}{\widetilde{\text{Cov}}(\mbf{q})}\left(\widetilde{\mathbb{E}}[\mbf{q}\, \mbf{q}_{\cdot+\tau}]-\widetilde{\mathbb{E}}[\mbf{q}]\, \widetilde{\mathbb{E}}[\mbf{q}_{\cdot+\tau}]\right).
\end{align}
where $\widetilde{\text{Cov}}(\mbf{q})$ denotes the covariance of $\mbf{q}$ across time as well as the ensembles; and $\widetilde{\mathbb{E}}$ denotes the time average with respect to the ensembles. Here, $\widetilde{\mathbb{E}}[\mbf{q}]$
approximates $\widehat{\mathbb{E}}[\mbf{q}]$, which in the continuous case is defined as
\[
\begin{split}
\widehat{\mathbb{E}}[\mbf{q}]:
=\lim_{T\to\infty}\frac{1}{M}\sum_{\ell =1}^{M}\left(\frac{1}{T}\int_0^T  \mbf{q}^{(\ell)}_t\, dt \right)
\quad \text{and}\quad
\widehat{\mathbb{E}}[\mbf{q}_{\cdot+\tau}]:
=\lim_{T\to\infty}\frac{1}{M}\sum_{\ell =1}^{M}\left(\frac{1}{T}\int_0^T  \mbf{q}^{(\ell)}_{t+\tau}\, dt \right).
\end{split}
\]
In the discrete case, an empirical approximation for uniform time-stepping size is defined as
\[
\widetilde{\mathbb{E}}[\mbf{q}]:=\frac{1}{M}\sum_{\ell=1}^{M} \frac{1}{N}\sum_{n=1}^{N} \mbf{q}_{t_n}^{(\ell)}
\quad \text{and}\quad
\widetilde{\mathbb{E}}[\mbf{q}_{\cdot+\tau}]:=\frac{1}{M}\sum_{\ell=1}^{M} \frac{1}{N}\sum_{n=1}^{N} \mbf{q}_{t_n+n_\tau}^{(\ell)},
\]
with $n_\tau$ being fixed and satisfying $\tau =\delta n_\tau.$

Based on the definition of $\widetilde{E}$, both  $\widetilde{\mathbb{E}}[\mbf{q}\,\mbf{q}_{\cdot+\tau}]$ and $\widetilde{\text{Cov}}$ are also computed via the empirical approximation across time and ensembles:
\[
\begin{split}
\widetilde{\mathbb{E}}[\mbf{q}\,\mbf{q}_{\cdot+\tau}]:=&\frac{1}{M}\sum_{\ell=1}^{M} \frac{1}{N}\sum_{n=1}^{N} \mbf{q}_{t_n}^{(\ell)} \mbf{q}_{t_n+n_\tau}^{(\ell)},
\quad\text{and}\quad
\widetilde{\text{Cov}}(\mbf{q}): =\widetilde{\mathbb{E}} \left[\left(\mbf{q}-\widetilde{\mathbb{E}}[\mbf{q}]\right) \left(\mbf{q}-\widetilde{\mathbb{E}}[\mbf{q}]\right) \right].
\end{split}
\]
\noindent In the following examples, we will consider two cases: 1) values of $\{\varepsilon_k\}$ are explicitly known, and 2) values of $\{\varepsilon_k\}$ are not explicitly known.

\subsection{Parameters for numerical tests}
Here, we list some numerical parameters which are commonly used in the following numerical tests. For 1D examples,  we fix the total number of time steps to be $N_t = 5\times 10^9$ and we set the scalar friction constant to be $\gamma =0.1$.  For 2D examples, we fix the total number of time steps to be $N_t = 2\times 10^7$, and use the Algorithm~\ref{Algorithm_summary} to match the friction matrix and diffusion coefficient according to \eqref{FDR}. In all examples of both dimensions, the time step size $\delta$ will be chosen according to the scale of interest $\varepsilon_k$.
When considering the effectiveness of surrogate Hamiltonian \eqref{surr_Hamiltonian}, we only generate one trajectory to compare the results of phase portrait. On the other hand, when comparing the Langevin system \eqref{Surr_Langevin}, we will use standard normal distributed initial conditions for $(\mbf{q},\,\mbf{p})$ to generate $M=4000$ independent trajectories on time interval $[0,\,T]$.

\subsection{One-dimensional examples with values of scales known}

\begin{enumerate}
    \item {\bf Test 1: quadratic potential.} We first consider a simple three-scale function, given by a quadratic $\mathcal{O}(1)$ component modulated by two smaller periodic scales
    \begin{equation}\label{test1}
    V(q)=\frac{q^2}{2}+\varepsilon_1\sin(q/\varepsilon_1)+\varepsilon_2\sin(q/\varepsilon_2) \text{ with } (\varepsilon_1, \, \varepsilon_2)=(0.05, \, 0.001).
    \end{equation}
    We run the scheme \eqref{verlet_scheme} with various time step sizes $\delta$ and plot the learnt functions in Figure~\ref{Fig:x2Potential_vs_step}. Notice that for $\delta=0.5$, the learned potential greatly resembles $U_0=\frac{q^2}{2}$, and for $\delta=0.065$, the learned potential greatly resembles $U_1=\frac{q^2}{2}+0.05\sin(q/0.05)$.
    Clearly, numerical step size as $\delta$ acts as a scaling filter to select the resolution of observed data.
    \begin{figure}[htp!]
    \centering
      \subfigure[Reference $\frac{q^2}{2}$]{\includegraphics[width=0.35\textwidth, height = 4cm ]{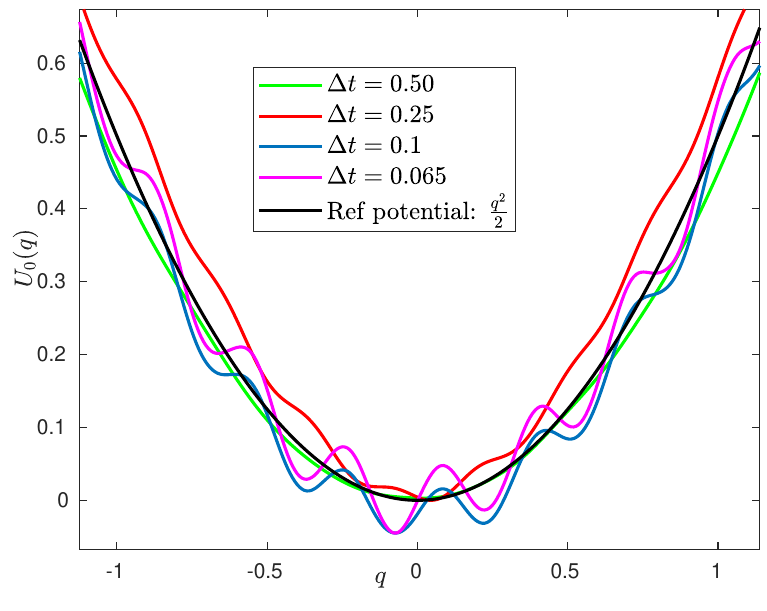}}
      \subfigure[Reference $\frac{q^2}{2}+\varepsilon_1\sin(q/\varepsilon_1)$]{\includegraphics[width=0.35\textwidth, height = 4cm ]{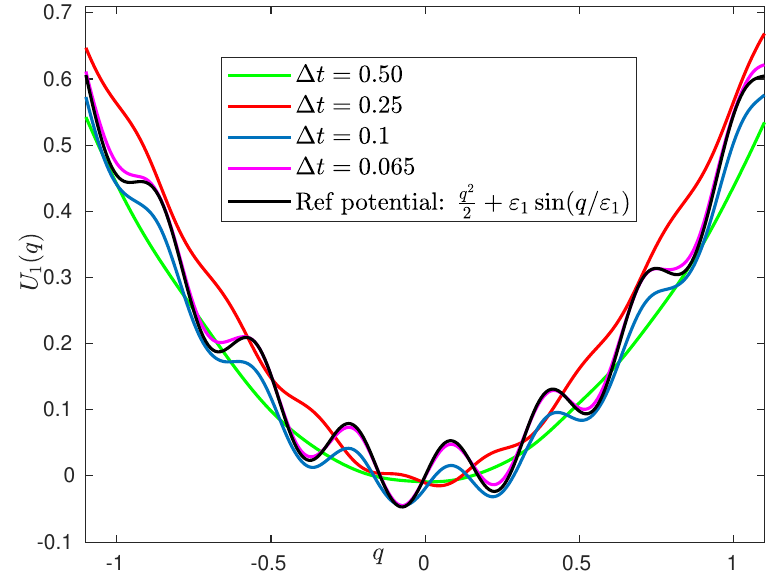}}
    \caption{The plots of learnt functions from simulated data with various time step sizes $\delta$ and the reference function ${U}_k(q)$. The exact function of full scales is $V(q)=\frac{q^2}{2}+\varepsilon_1\sin(q/\varepsilon_1)+\varepsilon_2\sin(q/\varepsilon_2)$. Parameters are $(\varepsilon_1, \, \varepsilon_2)=(0.05, \, 0.001)$.
    }
    \label{Fig:x2Potential_vs_step}
\end{figure}

\noindent
Based on learnt $U_0$, we run the surrogate Hamiltonian \eqref{surr_Hamiltonian} and the surrogate Langevin system \eqref{Surr_Langevin} using ${U}_0(q)$ with $\delta = 0.1$, and compare the simulations results using $V(q)$ with $h=5e-4$ for the Hamiltonian and Langevin dynamics, respectively.
As mentioned above, we compare the phase portrait for the Hamiltonian system, the equilibrium distribution, mean trajectory \eqref{eq:mean_traj} and normalized ACF \eqref{eq:ACF}. The results are summarized in Figure~\ref{Fig:x2_fit_vs_true}. For this quadratic example, we see that the surrogate models using $U_0(k)$ work well regarding capturing both statistics.

\begin{figure}[htp!]
    \centering
    \subfigure[Phase portrait]{ \includegraphics[width =0.4 \textwidth, height = 3.4 cm]{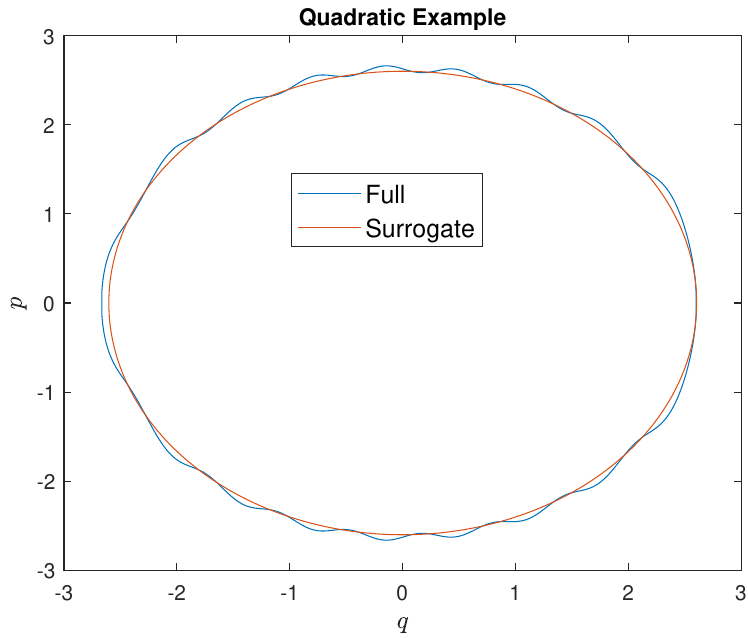}}
     \subfigure[Equilibrium dist]{ \includegraphics[width =0.4 \textwidth, height = 3.4 cm]{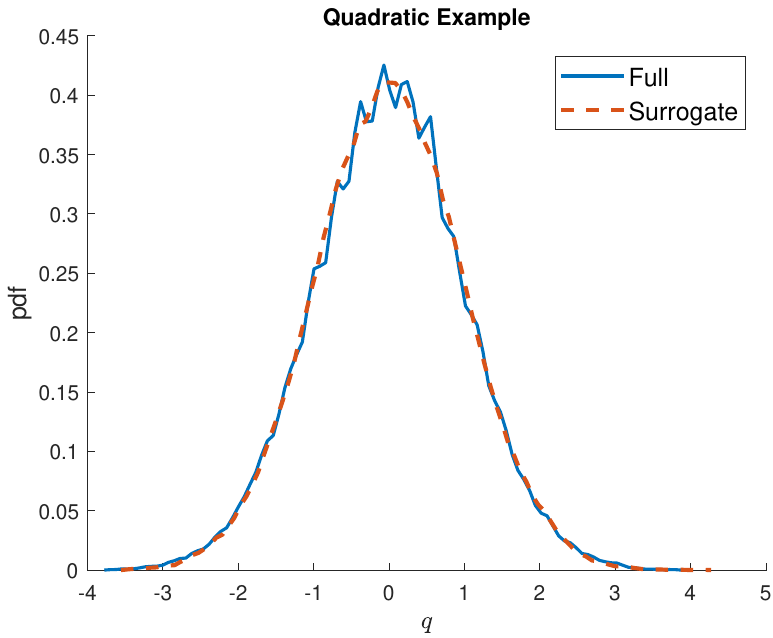}}
      \subfigure[Mean trajectory]{ \includegraphics[width =0.4 \textwidth, height = 3.4 cm]{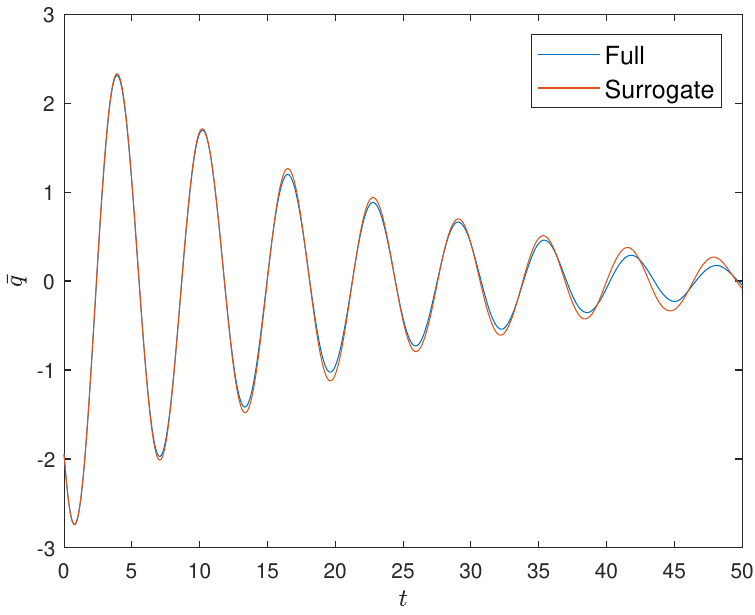}}
     \subfigure[Normalized ACF]{ \includegraphics[width =0.4 \textwidth, height = 3.4 cm]{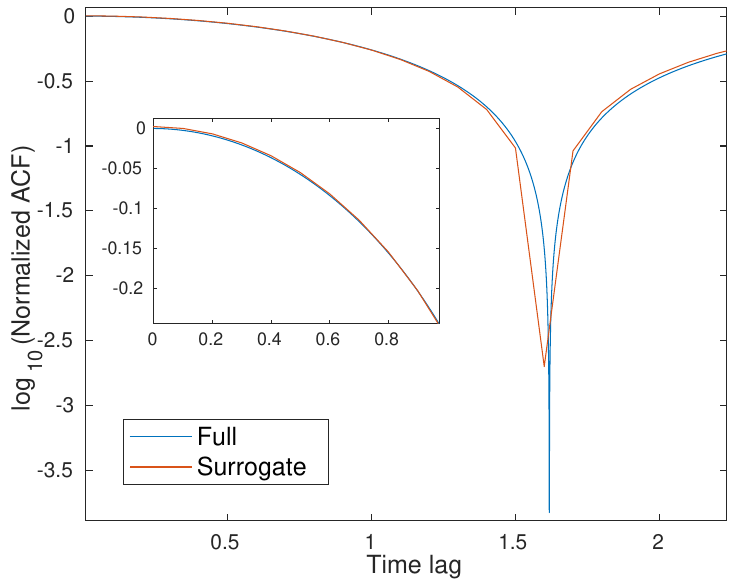}}
     \vspace{-0.2 cm}
    \caption{{\bf Fig (a):} Comparison on the surrogate Hamiltonian \eqref{surr_Hamiltonian} via $U_0(q)$ with $\delta = 0.1$ and via $V(q)$ with $h=5e-4$. {\bf Fig (b)-(d):} Comparison on the surrogate Langevin \eqref{Surr_Langevin} via $U_0(q)$ with $\delta=0.1$  and via $V(q)$ with $h = 5e-4$. $V(q)=\frac{q^2}{2}+\varepsilon_1\sin(q/\varepsilon_1)+\varepsilon_2\sin(q/\varepsilon_2)$. Parameters are $(\varepsilon_1, \, \varepsilon_2)=(0.05, \, 0.001)$ and $T=50$, initial distributions for $(q,\,p)\sim \mathcal{N}(-2,\,1)$.
    }
    \label{Fig:x2_fit_vs_true}
\end{figure}

\item {\bf Test 2: double-well potential.}
We next consider a three-scale function with its macroscopic component being a non-convex double-well function, instead of the previously considered convex quadratic function, and assume that we know the two small scales explicitly
\begin{equation}\label{test2}
V(q)=(q-1)^2\,(q+1)^2/4+\varepsilon_1\sin(q/\varepsilon_1)+\varepsilon_2\sin(q/\varepsilon_2) \text{ with } (\varepsilon_1, \, \varepsilon_2)=(0.025, 0.001).
\end{equation}
Similar to the quadratic example, we learnt $U_k(q)$ at each scale by choosing different time-stepping size $\delta$. The results of learnt $U_k$ with respect to various step sizes are plotted in Figure~\ref{Fig:DWPotential_vs_step}. We do capture scales of $\varepsilon_0$ and $\varepsilon_1$ by employing suitable $\delta$ for \eqref{verlet_scheme}.

    \begin{figure}[hpt!]
    \centering
      \subfigure[Reference $\frac{(q^2-1)^2}{4}$]{\includegraphics[width=0.35\textwidth, height = 4cm ]{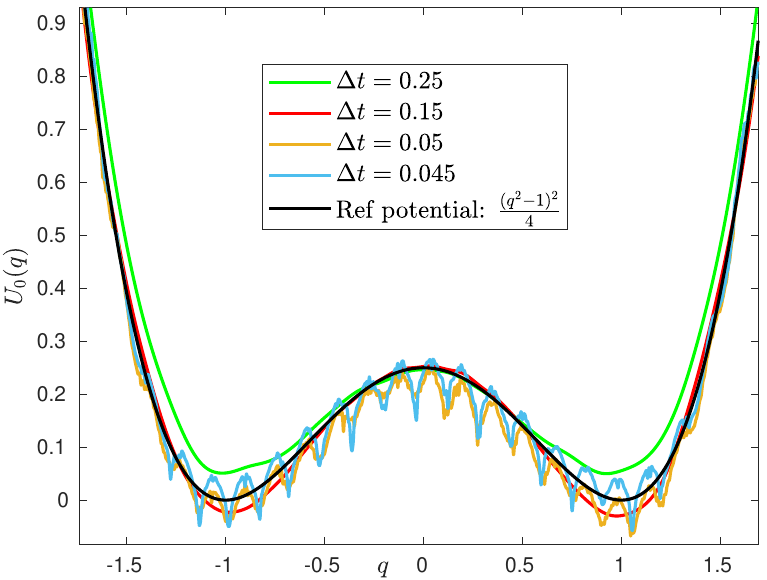}}
      \subfigure[Reference $\frac{(q^2-1)^2}{4}+\varepsilon_1\sin(q/\varepsilon_1)$]{\includegraphics[width=0.35\textwidth, height = 4cm ]{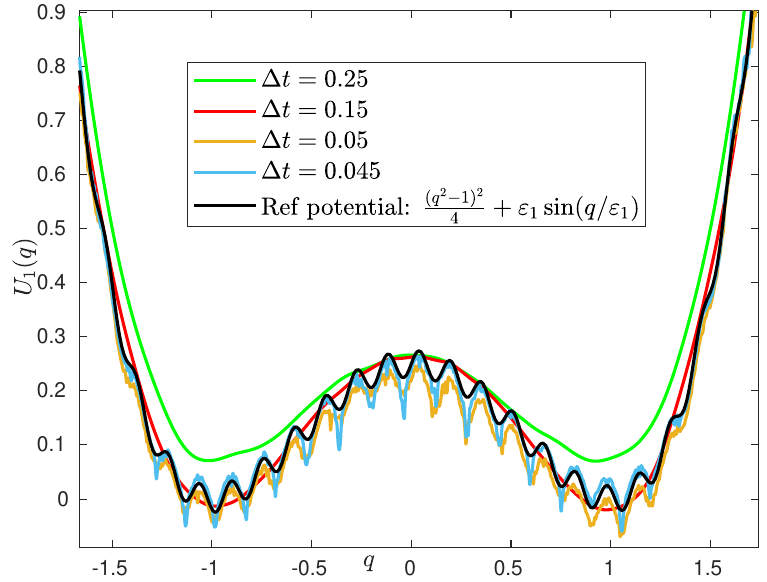}}
        \vspace{-0.2 cm}
    \caption{The plots of learnt functions from simulated data with various time step sizes $\delta$ and the reference function ${U}_k(q)$. The exact function of full scales is $V(q)=\frac{(q^2-1)^2}{4}+\varepsilon_1\sin(q/\varepsilon_1)+\varepsilon_2\sin(q/\varepsilon_2)$. Parameters are $(\varepsilon_1, \, \varepsilon_2)=(0.025, \, 0.001)$.     }
    \label{Fig:DWPotential_vs_step}
\end{figure}

We next compare the performances of  surrogate models via $U_0$ verse $V(q)$ for both Hamiltonian and Langevin. The results are summarized in Figure~\ref{Fig:DW_true_vs_fit}. We still observe very good match from the plots.

\begin{figure}[htp!]
       \centering
            \subfigure[Phase portrait] {\includegraphics[width=0.4\textwidth, height = 4 cm]{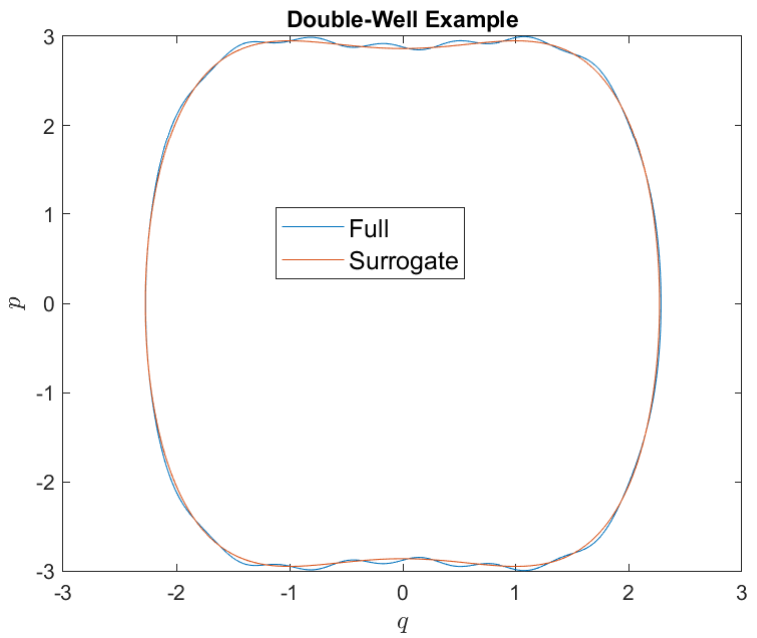}}
        \subfigure[Equilibrium dist]              {\includegraphics[width=0.4\textwidth, height = 4 cm]{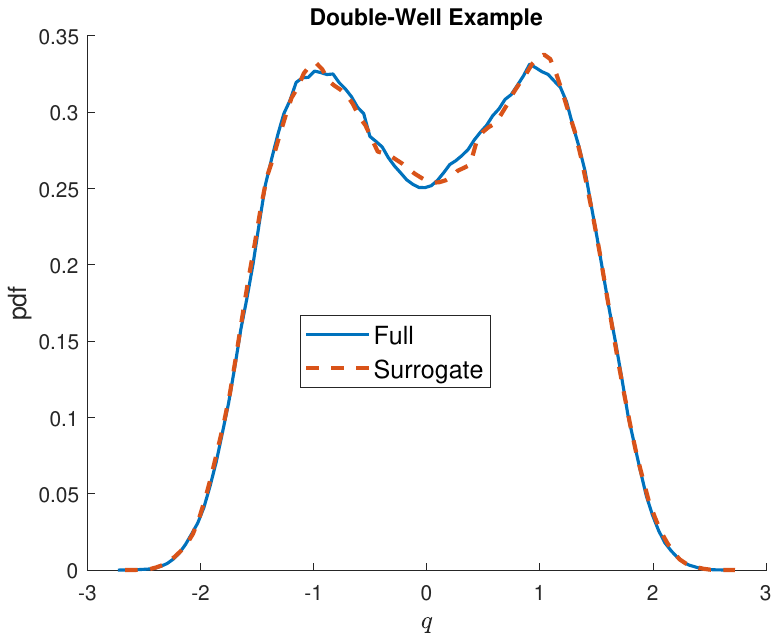}}
             \subfigure[Mean trajectory] {\includegraphics[width=0.4\textwidth, height =4 cm]{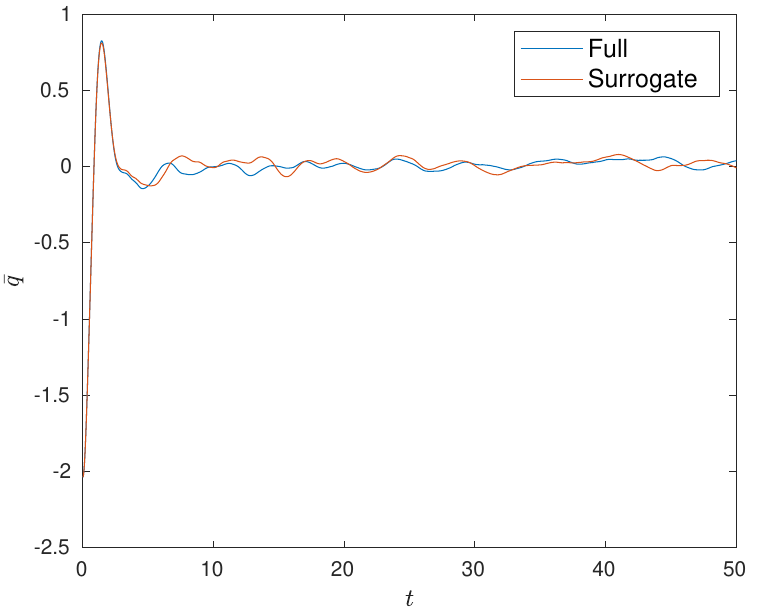}}
        \subfigure[Normalized ACF]              {\includegraphics[width=0.4\textwidth, height = 4 cm]{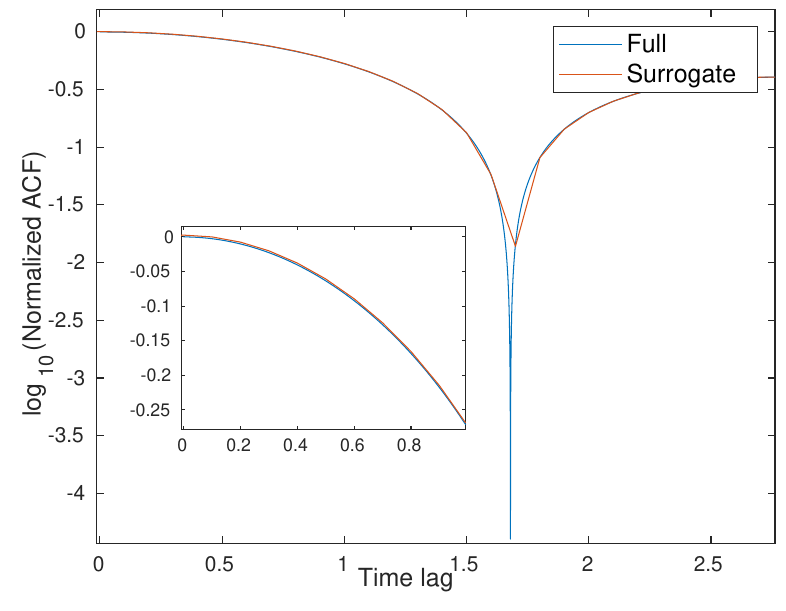}}
        \vspace{-0.2 cm}
               \caption{{\bf Fig (a):} Comparison on the surrogate Hamiltonian \eqref{surr_Hamiltonian} via $U_0(q)$ with $\delta = 0.1$ and via $V(q)$ with $h=5e-4$. {\bf Fig (b)-(d):} Comparison on the surrogate Langevin \eqref{Surr_Langevin} via $U_0(q)$ with $\delta = 0.1$ and via $V(q)$ with $h=5e-4$.  $V(q)=(q^2-1)^2/4+\varepsilon_1\sin(q/\varepsilon_1)+\varepsilon_2\sin(q/\varepsilon_2)$ with $ (\varepsilon_1, \, \varepsilon_2)=(0.025, 0.001)$.   $T=50$, initial distributions for $q\sim \mathcal{N}(-2,\,1)$ and $p\sim \mathcal{N}(0,\,1)$.
      }\label{Fig:DW_true_vs_fit}
   \end{figure}
\end{enumerate}

\subsection{1D example with lots of scales but no clear separation between adjacent scales}
\label{sec:noScaleSeparation}
Practical problems are oftentimes multiscale but without a sharp scale separation (two famous examples are fluid and molecular dynamics). Typically there is a wide range of scales, where adjacent two scales are not well separated but the largest and smallest scale are very well separated. To consider a simplified version where the focus is just to decompose a multiscale function, we investigate again problem \eqref{multiscale_V}
, however with
$1 \gtrsim \varepsilon_1 \gtrsim \cdots \gtrsim \varepsilon_K$ and $1\gg \varepsilon_K$. Suppose $\delta$ is chosen at scale $k$, i.e. $\delta\approx \varepsilon_k$, the effective function should not be exactly $V_0(q)+\sum_{j=1}^k V_{j,\varepsilon_j}$, but with some additional correction as percolated from nearby smaller scales $\varepsilon_{k+1}$, etc. However, there is no clear or unique theoretical justification on how this correction works. We would like to inspect what our approach would produce in this case, and how the corresponding surrogate model approximates the original full model.

\begin{figure}[htp!]
    \centering
    \includegraphics[width = 0.4\textwidth, height = 4 cm]{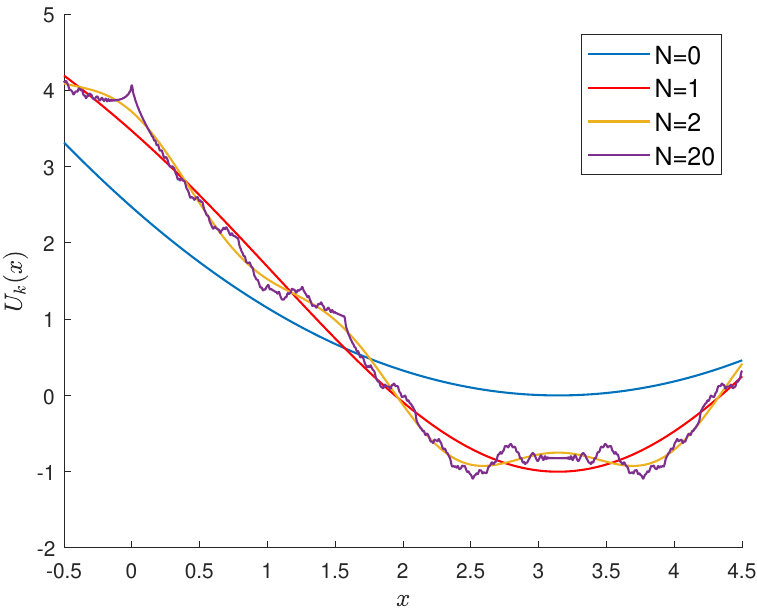}
    \vspace{-.2 cm}
    \caption{Zoomed-in plot of $V(q)=\frac{(q-\pi/2)^2}{4}+\sum_{i=1}^{N}\frac{\cos(i^2\, q)}{i^2}$ with different values of $N$.}
    \label{fig:cos_demo}
\end{figure}

\noindent
For this purpose, we consider a toy test problem
\begin{equation}\label{cos_pot}
V(q)=\frac{(q-\pi/2)^2}{4}+\sum_{i=1}^{N}\frac{\cos(i^2\times q)}{i^2}.
\end{equation}
A demonstration of the exact potential at different scales can be found in Figure~\ref{fig:cos_demo}. Notice that different $i$'s (besides the 1st term) only nominally correspond to different scales. While adjacent scales are close to each other, when $N$ is large, $V(\cdot)$ clearly exhibits a wide range of well separated scales. In addition, one cannot directly read off $\{\varepsilon_k\}$ and $U_k(q)$ from the expression of $V(q)$, because different $i$'s can contribute to the same scale (and other scales too) due to a lack of clear scale separation. Instead, it is reasonable to postulate the existence of some hidden functions that represent different scales.

\begin{figure}[htp!]
       \centering
       \subfigure[$\delta = 0.6500$] {\includegraphics[width=0.32\textwidth, height = 3.4 cm]{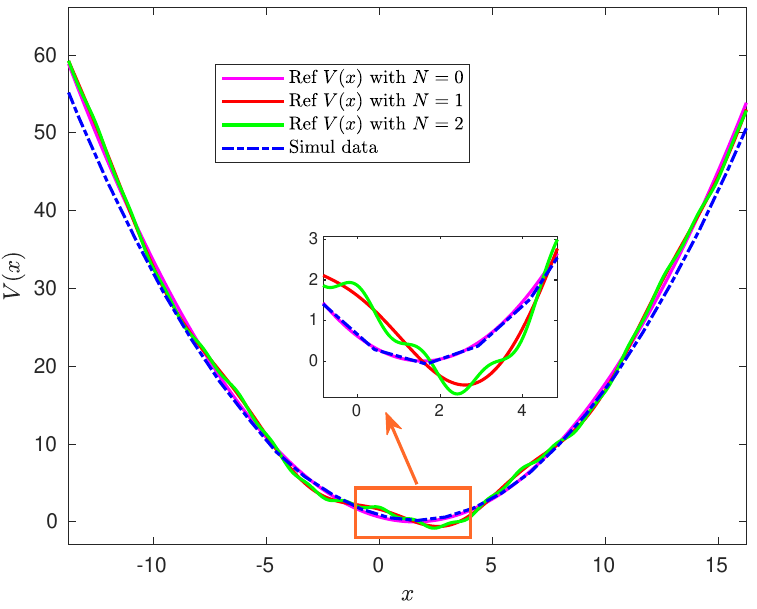}}
       \subfigure[$\delta = 0.2601$] {\includegraphics[width=0.32\textwidth, height = 3.4 cm]{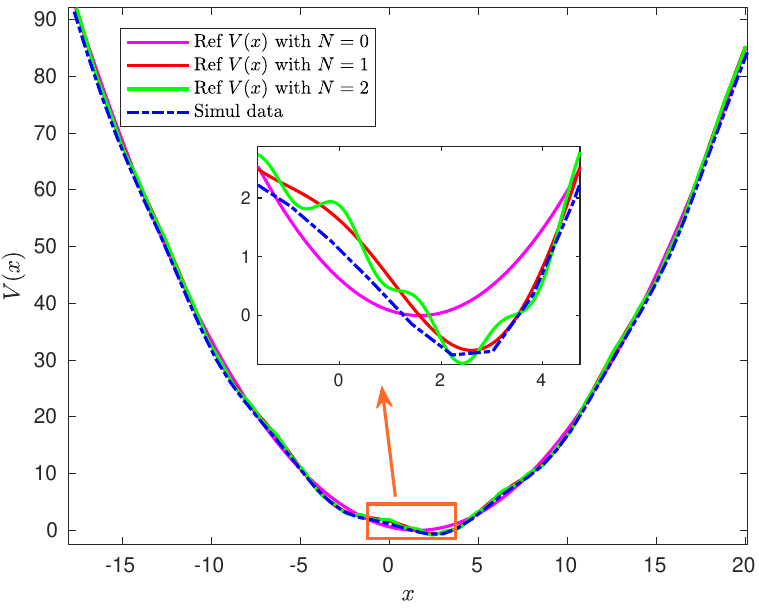}}
        \subfigure[$\delta = 0.1309$] {\includegraphics[width=0.32\textwidth, height = 3.4 cm]{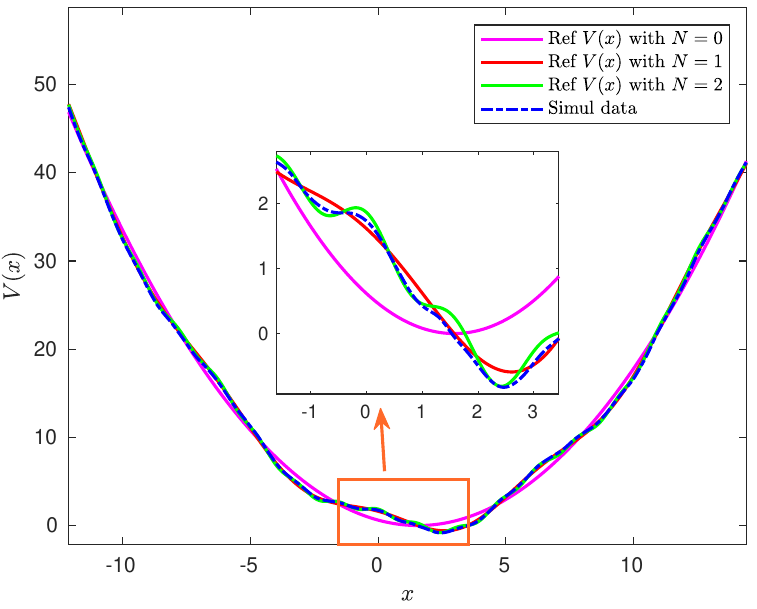}}
         \caption{Learning potentials $U_k(q)$ at different scales. The exact potential is $V(q):=\frac{(q-\pi/2)^2}{4}+\sum_{i=1}^{20}\cos(i^2\times q)/i^2$.
      }\label{Fig:cosHighF_Potlearnt}
\end{figure}

\begin{figure}[htp!]
       \centering
            \subfigure[Phase portrait] {\includegraphics[width=0.4\textwidth, height = 4 cm]{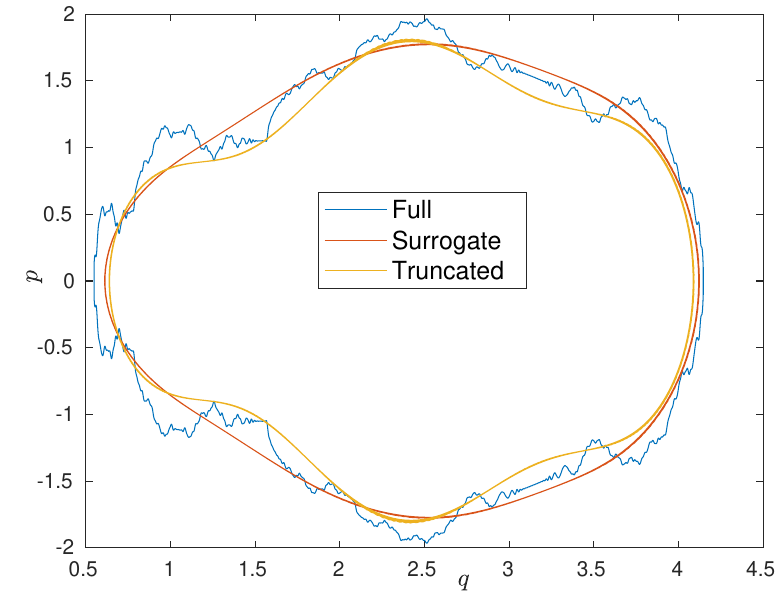}}
        \subfigure[Equilibrium dist $\pi_\delta(q)$]              {\includegraphics[width=0.4\textwidth, height = 4 cm]{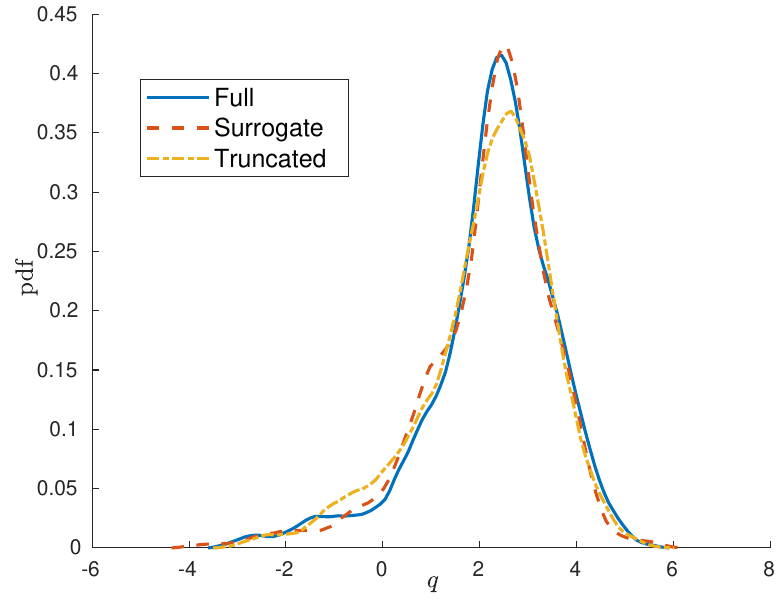}}
             \subfigure[Mean trajectory             ] {\includegraphics[width=0.4\textwidth, height =4 cm]{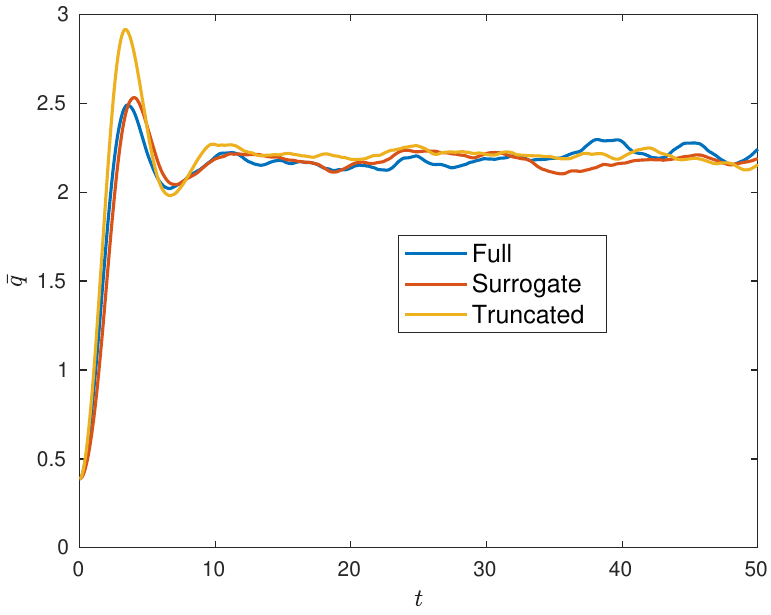}}
        \subfigure[Normalized ACF]              {\includegraphics[width=0.4\textwidth, height = 4 cm]{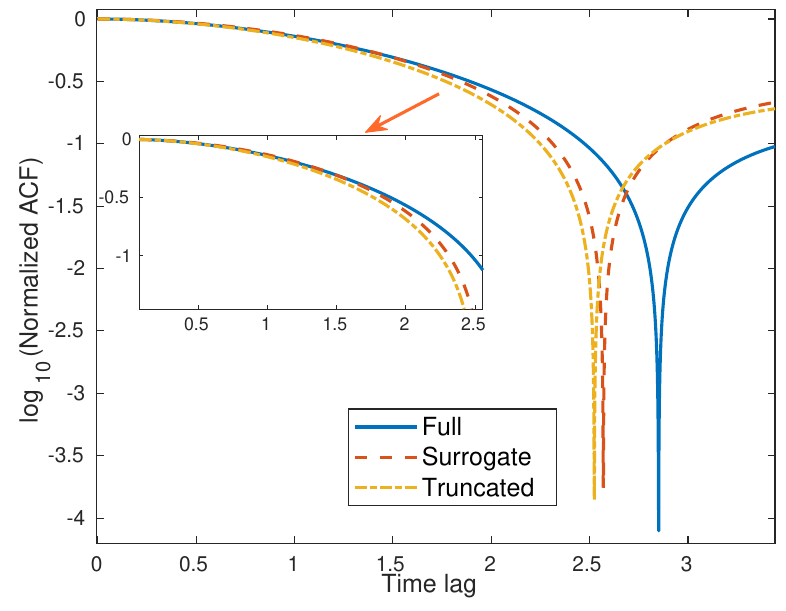}}
        \vspace{-0.2 cm}
               \caption{{\bf Fig (a):} Comparison on the surrogate Hamiltonian \eqref{surr_Hamiltonian} via $U_0(q)$ with $\delta = 0.1309$ and via $V(q)$ with $h=1e-3$. {\bf Fig (b)-(d):} Comparison on the surrogate Langevin \eqref{Surr_Langevin} via $U_0(q)$ with  $\delta = 0.1309$
               and via $V(q)$ \eqref{cos_pot} for $N=20$ as well as the truncated version $N=2$  with $h = 1e-3$.  The initial distributions for $q\sim \mathcal{U}(-1/2,\,1/2)+0.38$ and $p\sim \mathcal{U}(-1/2,\,1/2)$.
      }\label{Fig:MultiCos_true_vs_fit}
   \end{figure}

\noindent
To learn the mesocale $\{U_k(q)\}$ at various scales, we fix $\gamma=0.5$, mass $M=1$ and choose the total number of time steps to be $N_t = 5\times 10^9$. We find several scales of potential functions via tuning the step sizes $\delta$. The results are present in Figure~\ref{Fig:cosHighF_Potlearnt}. We note from simulation that there is a macroscopic scale $\varepsilon_0$ around $\varepsilon_0\sim \delta=0.65$ as the learnt function is close to the referenced $N=0$. The two following scales corresponding to $N=1$ and $N=2$ are around $\varepsilon_1 \sim \delta = 0.2601$ and $\varepsilon_2 \sim \delta = 0.1309$. The numerical separation between $\varepsilon_1$ and $\varepsilon_2$ is not as clear as those at $\varepsilon_0$ (i.e., $N=0$). In fact, $\varepsilon_1$ and $\varepsilon_2$ are just weakly separated due to the narrow value gap in $V(q)$ between $N=1$ and $N=2$. Moreover, as demonstrated in Figure~\ref{Fig:MultiCos_true_vs_fit}, the accuracy of surrogate Hamiltonian and Langevin models constructed by our learning algorithm is better than those constructed by a simple truncation of objective function due to the mixed scales among the truncation.

\subsection{2D quadratic function with anisotropic small scale}
We now consider multi-dimensional examples. Throughout these examples, we will assume zero knowledge about how $V$ decomposes in scales even though we wrote down its expression analytically.

\begin{figure}[htp!]
    \centering
      \subfigure[Ref $V_0$ and learning results]{\includegraphics[width=0.4\textwidth, height = 4cm ]{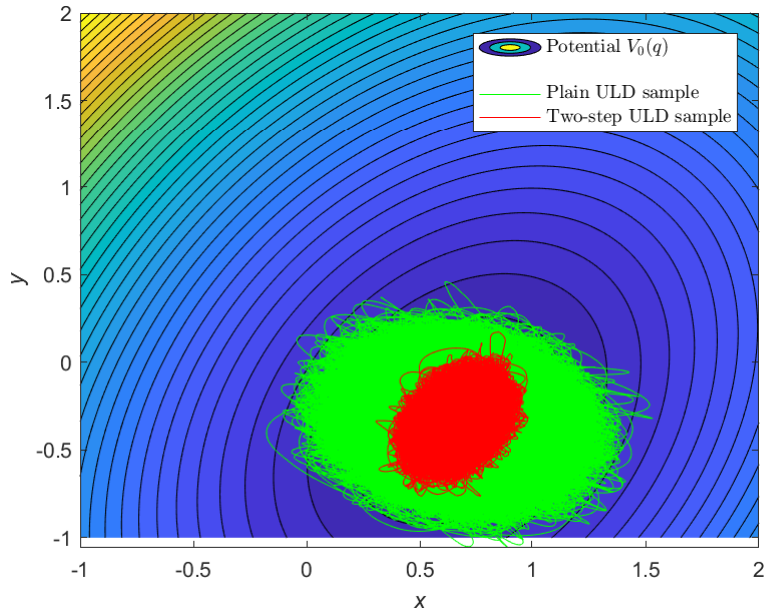}}
      \subfigure[Zoomed-in view]{\includegraphics[width=0.4\textwidth, height = 4cm ]{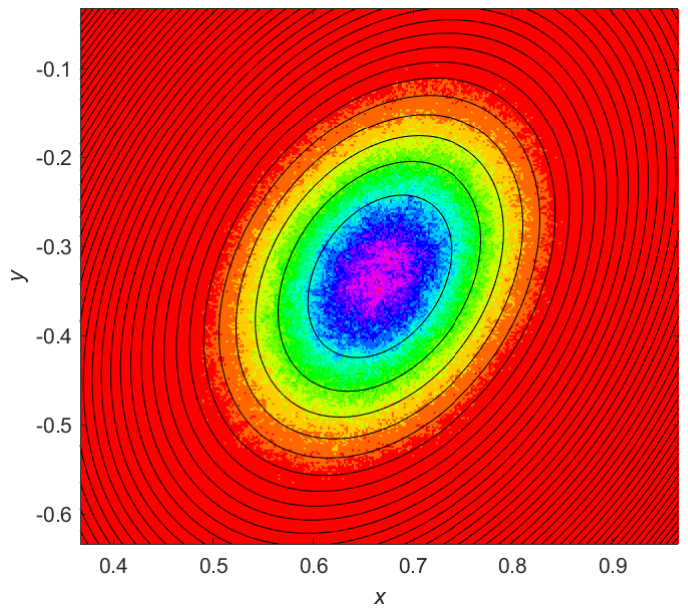}}
      \caption{Learning results of macroscopic $U_0$. The exact potential is given in \eqref{2d_quadratic} and the step size is $\delta = 0.05$. {\bf Left:} the green data is generated using scalar friction $\gamma I $ only; and the red data is generated by the two-stage Algorithm~\ref{Algorithm_summary}. {\bf Right:} zoomed-in view of the two-stage learning results. In the legend, ULD denotes the under damping-loaded dynamics \eqref{Newton_sys}.
      }\label{Fig:Ellipse}
   \end{figure}

The first one is a two-dimensional quadratic potential, $\mbf{q}=(x,\,y)$ with one known anisotropic small scale
\begin{equation}\label{2d_quadratic}
\begin{aligned}
 V(x,y) = &V_0+V_1\\
 =& \frac{1}{4}(2x+y-1)^2+(x-y-1)^2  +  \varepsilon\bigg(\sin(x/\varepsilon)+\sin\left((x+y)/\varepsilon\right)\bigg),
 \quad \varepsilon=  10^{-5}.
 \end{aligned}
\end{equation}
We apply the two-stage Algorithm~\ref{Algorithm_summary} with $\delta = 0.05$ and $N_t= 2\times 10^7$ to learn the macroscopic potential $U_0=V_0(x,y)$. During stage 1, the scalar friction $\gamma$ is chosen to be $\gamma=0.1$. During stage 2, we estimate the variance of small scale contributions and set $\Gamma$ accordingly. We then compare the learning function from using $\gamma I$ only with that from the two-stage algorithm. The results are summarized in Figure~\ref{Fig:Ellipse}. We emphasize the mismatch of green data is generated without normalizing the variance of small scale,
whereas the data from two-stage simulation can capture the macroscopic $U_0$ very well.

   \subsection{2D M{\"{u}}ller-Brown potential with one known anisotropic small scale}
   M{\"{u}}ller-Brown potential is a common toy model used in molecular dynamics research, due to its high nonlinearity and the existence of multiple minima.
   We now consider the macrosopic potential to be a modification of M{\"{u}}ller-Brown potential, and the microscopic scale  to be some anisotropic toy function. More precisely, consider
   \begin{equation}\label{Modified_Muller}
       \begin{aligned}
       V(x,y)=&V_0(x,y)+V_1(x,y),\\
       V_0(x,y)=& 0.1\times\bigg(V_q(x,y)+ V_m(x,y)\bigg),\\
         V_1(x,y)=& \varepsilon\bigg(\sin\big(x/\varepsilon\big)+\sin\big((-x+y)/\varepsilon\big)\bigg),\quad \varepsilon=10^{-5},
       \end{aligned}
   \end{equation}
   where  $V_q(x,y)$ is a simple quadratic potential
   $V_q(x,y)=35.0136(x-x_c^*)^2+ 59.8399 (y-y_c^*)^2$, and
   $V_m(x,y)$ denotes the M{\"{u}}ller-Brown potential
   \[
   \begin{split}
   V_m(x,y)=&-200\exp\big(-(x-1)^2-10\,y^2\big)    -100\exp\big(-x^2-10(y-0.5)^2\big)\\
   &\quad -170\exp\big(-6.5(x+0.5)^2+11(x+0.5)(y-1.5)-6.5(y-1.5)^2\big)\\
   &\quad\quad +15\exp\big(0.7(x+1)^2+0.6(x+1)(y-1)+0.7(y-1)^2\big).
   \end{split}
   \]
Here, $(x_c^*,\,y_c^*)$ denotes the center of the middle well of $V_m$. Notice that $V_q(x,y)$ is introduced so that the depths of all three wells of the modified  M{\"{u}}ller-Brown (MB) function are better leveled to each order. This is because, the original M{\"{u}}ller-Brown without modification has three local minimina with very different function values, and therefore at the thermal equilibrium (i.e. Gibbs distribution) some potential well is exponentially less likely visited than the others, which both reduces the learning accuracy of the corresponding mode and makes visualization difficult.

\begin{figure}[htp!]
    \centering
      \subfigure[Ref $V_0$ and learning results  ]{\includegraphics[width=0.44\textwidth, height = 4cm ]{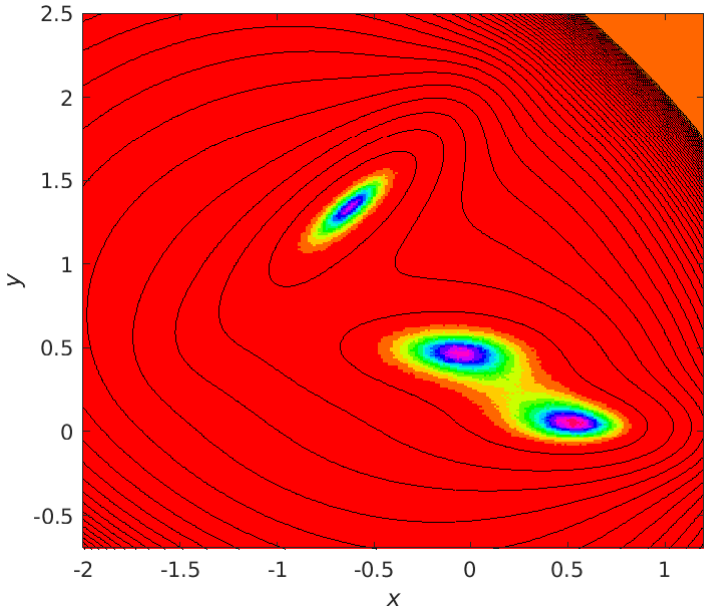}}
      \subfigure[Zoomed-in view]{\includegraphics[width=0.44\textwidth, height = 4cm ]{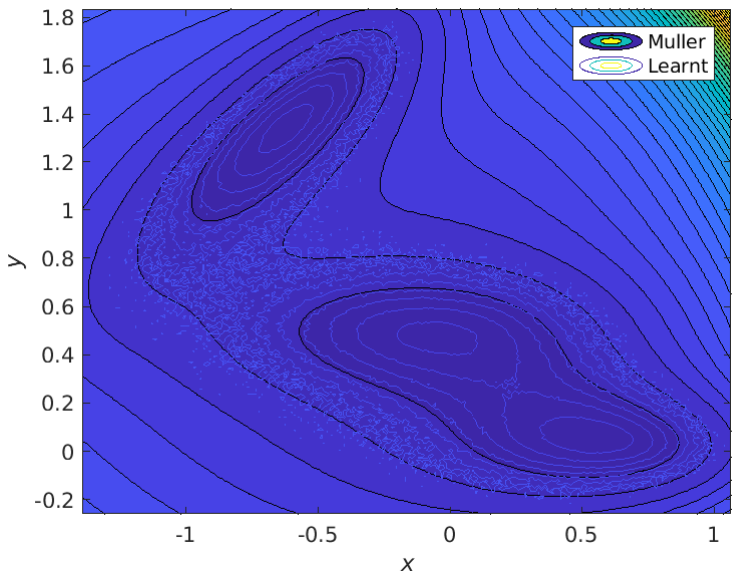}}
      \caption{Learning results of macroscopic $U_0$ at three wells. The exact potential is given in \eqref{Modified_Muller} and the simulation step size is $\delta = 0.05163$. The three wells of $U_0(x,y)$ are learnt simultaneously.
      }\label{Fig:Muller}
   \end{figure}

We apply the two-stage Algorithm with $\delta = 0.05163$
and $N_t= 1\times 10^8$ to learn the macroscopic scale $U_0=V_0(x,y)$. In the stage 1, the scalar friction $\gamma$ is chosen to be $\gamma=0.014$.
The results are summarized in Figure~\ref{Fig:Muller}. Note that we can learn all three energy wells at the same time with the existence of anisotropic small scale components.

\section{Conclusion}\label{sec:conclusion}
In this study, we introduced a novel algorithm for decomposing and learning every scale of a given multiscale objective function across multiple dimensions. Our method leverages the controllable algorithmic implicit bias inherent in the stochasticity of large learning rate gradient descent,
thus facilitating the automatic generation of distinct data sets across varying scales. Upon  successfully decomposing the multiscale objective function, many applications can be enabled. As a demonstration, we constructed a surrogate model aligned with the equilibrium distribution and dynamic mean path and auto-correlation of data produced by the original objective function at the scale of interest.

Additionally, we devised a two-stage algorithm to segregate variables of interest within an anisotropic stochastic forcing environment. By adjusting the system's mass tensor, we estimated the covariance of the stochastic forcing, and an appropriate friction matrix could be chosen to bring the equilibrium distribution back towards a Gibbs distribution, which possesses an analytical expression.

Meanwhile, the proposed strategy has certain limitations. The first limitation is the statistical accuracy due to finite samples. Detailed choices for the function fitting also matter; for example, the effective potential can only be resolved to the scale dictated by the bin width, and which (parameterized) model to use for the regression makes a difference. Although these are all generic and standard problems, we would still like to point them out, because a consequence is, the current implementation requires a large amount of (self-generated) data. The second limitation is the necessity of introducing hyper parameters when generating data, such as the friction constant $\gamma\in (0,\, 0.1]$ and proper range of $\hat{\beta}\in(0.2,\,1.25)$ in the 1D simulation. We do observe that the learning results tend to better when we use that range of parameters.

Possible future directions include 1) to explore more statistical tools to have better estimation of $\pi_\delta$ and more robust regression of $U_k(\mbf{q})$, and 2) to derive \textit{a priori} error estimate of learning $U_k$ in terms of various problem and hyper parameters settings.

\section*{Acknowledgements}
X. Li’s is grateful for partial support by the NSF Award DMS-1847770 and the 2023 UNC Charlotte faculty research grant. M. Tao is grateful for partial support by the NSF Award DMS-1847802, NSF Award ECCS-1936776, Cullen-Peck Scholar Award, and Emory-GT AI. Humanity Award.

\bibliographystyle{plain}
{\small
\bibliography{Ref_ScaleDecomp}
}

\end{document}